\documentclass[reqno, 11pt]{amsart}

\usepackage{tikz}
\usepackage{pgfplots}
\usetikzlibrary{arrows,shapes,positioning,patterns}
\usetikzlibrary{decorations.markings}
\tikzstyle arrowstyle=[scale=1]
\tikzstyle directed=[postaction={decorate,decoration={markings,
		mark=at position .65 with {\arrow[arrowstyle]{stealth}}}}]
\tikzstyle reverse directed=[postaction={decorate,decoration={markings,
		mark=at position .65 with {\arrowreversed[arrowstyle]{stealth};}}}]

\usepackage[english]{babel}
\usepackage{amsmath}
\usepackage{amsthm}
\usepackage{amsfonts}
\usepackage{amssymb}
\usepackage{anysize}
\usepackage{hyperref}
\usepackage{appendix}
\usepackage{mathtools}
\usepackage{graphicx}
\usepackage{xcolor}
\usepackage{pdfpages}
\usepackage{ulem}

\newtheorem{definition}{Definition}[section]
\newtheorem{lemma}[definition]{Lemma}
\newtheorem{theorem}[definition]{Theorem}
\newtheorem{corollary}[definition]{Corollary}
\newtheorem{proposition}[definition]{Proposition}

\newtheorem{remark}[definition]{Remark}

	\title[Keller--Segel type with logarithmic sensitivity]{Kinks and solitons in linear and nonlinear--diffusion Keller--Segel type models with logarithmic sensitivity}

	\author[J. Campos]{Juan Campos}
	\address{Departamento de Matem\'atica Aplicada and  Research Unit ``Modeling Nature'' (MNat), Facultad de Ciencias. Universidad de Granada\\ 18071-Granada, Spain}
	\email{campos@ugr.es}
	\author[C. Garc\'ia]{Claudia Garc\'ia}
	\address{Departamento de Matem\'atica Aplicada and  Research Unit ``Modeling Nature'' (MNat), Facultad de Ciencias. Universidad de Granada\\ 18071-Granada, Spain. 
	}
	\email{claudiagarcia@ugr.es}
	\author[C. Pulido]{Carlos Pulido}
	\address{Departamento de Matem\'atica Aplicada and  Research Unit ``Modeling Nature'' (MNat), Facultad de Ciencias. Universidad de Granada\\ 18071-Granada, Spain}
	\email{cpulidog@ugr.es}
	\author[J. Soler]{Juan Soler}
	\address{Departamento de Matem\'atica Aplicada and  Research Unit ``Modeling Nature'' (MNat), Facultad de Ciencias. Universidad de Granada\\ 18071-Granada, Spain}
	\email{jsoler@ugr.es}
	{\thanks{This work has been partially supported by the State Research Agency of the Spanish Ministry of Science and FEDER-EU, project PID2022-137228OB-I00 (MICIU/AEI /10.13039/501100011033); by Modeling Nature Research Unit, Grant QUAL21-011 funded by Consejer\'ia de Universidad, Investigaci\'on e Innovaci\'on (Junta de Andaluc\'ia); and by Grant C-EXP-265-UGR23 funded by
Consejer\'ia de Universidad, Investigaci\'on e Innovaci\'on \& ERDF/EU
Andalusia Program. C.G. has been also supported by RYC2022-035967-I (MCIU/AEI/10.13039/501100011033 and
FSE+), and partially by Grant PID2022-140494NA-I00 and a 2024
Leonardo Grant LEO24-1-13722 for Scientific Research and Cultural Creation from the BBVA
Foundation.}}
	
	\subjclass[2010]{35A01, 335A15, 5B40, 35B44, 35C07, 35K57, 35K65, 35Q92, 92C17}
	\keywords{Traveling waves, Keller-Segel, Flux-saturated, Patterns in nonlinear parabolic systems, Solitons, Logarithmic sensitivity, Cross-diffusion}
	
\begin{document}
		
		\date{}
		
		\begin{abstract}
This paper investigates the existence of traveling--wave--type patterns in the Keller--Segel model with logarithmic sensitivity. We consider both the linear diffusion case and the nonlinear, flux-saturated diffusion of relativistic heat--equation type, providing a detailed comparison between the two regimes. Particular attention is devoted to traveling waves exhibiting compact support or support restricted to a half-line. We rigorously establish the existence of such patterns and highlight the qualitative differences arising from the choice of diffusion mechanism.
\end{abstract}
		
		\maketitle
		
		\tableofcontents

		\section{Introduction}
		The aim of this paper is to study traveling wave patterns for the Keller--Segel model with logarithmic sensitivity, both in the case of linear diffusion and in the case of flux-saturated nonlinear diffusion. Our analysis will focus in particular on the relativistic heat equation, a prototypical example of flux limitation. Our objective is to prove the existence of soliton-type traveling waves with compact support in both the linear and nonlinear diffusion regimes.

Chemotaxis refers to the directed motion of a biological species along gradients of a chemical signal. Classical examples include the formation and propagation of traveling bacterial bands toward oxygen \cite{Adler1, Adler2}, or the concentric ring waves observed in {\it E.~coli} colonies \cite{Brenner, Budrene}. The prototypical mathematical description of chemotaxis was introduced by Keller and Segel \cite{[KS71]}. In its general form, the model reads
		\begin{equation}\label{system-0} \left\{
			\begin{array}{ll}
				\displaystyle \partial_t u(t,x)=\partial_x \left\{u(t,x)\Phi\left(\frac{\partial_x u(t,x)}{u(t,x)} \right)-a u(t,x) \partial_x f(S)\right\},& x\in\mathbb{R}, t>0,\\ {} \\
				\delta \partial_t S(t,x)=\gamma \partial^2_{xx}S(t,x)+k(u,S), & x\in\mathbb{R}, t>0,\\ {}\\
				u(0,x)=u_0(x), &x\in\mathbb{R}.
			\end{array}\right.
		\end{equation}
		Here $u=u(t,x)$ denotes the cell density, and $S=S(t,x)$ the chemoattractant concentration. The parameter $a\ge0$ is the chemotactic sensitivity, $\gamma\ge0$ is the diffusion coefficient of the chemical, and $\delta\ge0$ describes the time scale for the chemical dynamics. In the classical Keller--Segel system one sets $\Phi(s)=s$, yielding linear diffusion, while the function $f$ encodes the chemosensitivity mechanism, and $k(u,S)$ represents chemical production and degradation.

Different choices of the chemosensitivity function appear in the literature: the linear law $f(S)=S$, the logarithmic law $f(S)=\log S$, and the receptor law $f(S)=S^m/(1+S^m)$ for $m\in\mathbb{N}$. The linear law with $k(u,S)=S-u$ leads to the {\it minimal} chemotaxis model \cite{Childress, Horstmann}. The logarithmic law, motivated by the Weber--Fechner principle, has been widely studied (see, e.g., \cite{Alt-Lauffenburger, Balding-McElwain, Dahlquist-Lovely-Koshland, [KS71]}). Although originally developed for chemotaxis, Keller--Segel type models now appear in diverse contexts including population dynamics, polymer science, and quantum cross-diffusion. We refer to \cite{Wang} for a survey on the logarithmic case.

In this work we focus on the logarithmic sensitivity $f(S)=\log S$ together with the choice $k(u,S)=u-\lambda S$ for $\lambda\ge0$. In this case system \eqref{system-0} becomes
		\begin{equation}\label{system} \left\{
			\begin{array}{ll}
				\displaystyle \partial_t u(t,x)=\partial_x \left\{u(t,x)\Phi\left(\frac{\partial_x u(t,x)}{u(t,x)} \right)-a\frac{\partial_x S(t,x)}{S(t,x)} u(t,x)\right\},& x\in\mathbb{R}, t>0,\\ {} \\
				\delta \partial_t S(t,x)=\gamma \partial^2_{xx}S(t,x)-\lambda S(t,x)+u(t,x), & x\in\mathbb{R}, t>0,\\ {}\\
				u(0,x)=u_0(x), &x\in\mathbb{R}.
			\end{array}\right.
		\end{equation}

		We will assume that $\Phi$ verifies
		\begin{enumerate}
			\item[{\bf(H1)}]  $
			\Phi\in C^2(\mathbb{R}), \, \Phi(-s)=-\Phi(s), \, \Phi'(s)>0, \, \forall s\in \mathbb{R}.
			$ 
		\end{enumerate}
		We assume that
\begin{enumerate}
\item[{\bf(H1)}] $\Phi\in C^2(\mathbb{R}),\quad \Phi(-s)=-\Phi(s),\quad \Phi'(s)>0\ \text{for all }s\in\mathbb{R}.$
\end{enumerate}
The classical Keller--Segel model corresponds to $\Phi(s)=\mu s$ for some viscosity coefficient $\mu>0$.  

We will compare the linear diffusion case with the flux-limited case. If
\[
\lim_{s\to\infty} \Phi(s)=\infty,
\]
as in the classical diffusion setting, no qualitatively new phenomena appear, apart from certain technical difficulties. In contrast, if
\begin{enumerate}
\item[{\bf(H2)}] $\displaystyle \lim_{s\to\infty} \Phi(s)<\infty,$
\end{enumerate}
the diffusion becomes flux-saturated, and new qualitative behaviours emerge.

A prominent example is the {\it relativistic heat equation}, corresponding to
		$$
		\Phi(s)=\mu \frac{s}{\sqrt{1+\left(\frac \mu c\right)^2 s^2}},
		$$
where $\mu,c>0$. In this case the limit in {\bf(H2)} equals $c$, which governs the maximal propagation speed of the support. We refer to {\bf(H)} to denote hypotheses {\bf(H1)}--{\bf(H2)}. The Larson-type operators
\[
\Phi(s)=\mu\,\frac{s}{\sqrt[p]{1+\left(\frac{\mu}{c}\right)^{p}|s|^p}},\qquad p\in(1,\infty),
\]
and the Wilson operator (the case $p=\infty$),
\[
\Phi(s)=\mu\,\frac{s}{1+\left(\frac{\mu}{c}\right)|s|},
\]
can also be treated by similar techniques. We refer to \cite{CCCSS-SUR} for an extensive survey on flux-saturated diffusions. These operators typically satisfy a sublinear growth condition
\[
|\Phi(s)|\le a|s|+b,\qquad a,b>0.
\]

Models with flux limitation have been studied from many perspectives: well-posedness and entropy solutions \cite{ACMSV, ARMAsupp, ACM2005, CCCSS-EMS}, waiting-time phenomena \cite{CCCSS-EMS, Giacomelli}, hydrodynamic limits for multicellular systems \cite{[BBNS07],[BBNS10],[BBNS10B], BBTW, BW1, BW2}, connections with FKPP or porous medium equations \cite{Campos2016, preprintVicent}, and applications ranging from astrophysics \cite{Levermore81, MM} to morphogenesis, sound propagation, and tumour growth \cite{Blanco, CMSV, Campos, CKR, Conte, VGRaS}.
		
		In this work we study traveling wave solutions of the form
\begin{equation}\label{def-tw}
u(t,x)=u(x-\sigma t),\qquad S(t,x)=S(x-\sigma t),
\end{equation}
with speed $\sigma>0$. Traveling waves are fundamental in understanding the emergence of propagating patterns, and have been studied extensively for Keller--Segel type models (see \cite{[HP],[KS71],[K80],[PA], sanchez} and the references therein).

We show that the flux-saturated Keller--Segel system exhibits behaviours absent from the classical model. Figure~\ref{types} illustrates the possible wave profiles in both cases. In the flux-saturated regime the traveling waves may exhibit compact support with sharp edges and infinite slopes at the boundaries, a phenomenon impossible in the classical diffusive case.
		\begin{figure}[h] \label{types}
			\begin{center}
				\begin{tikzpicture}[scale=1.3]
					\draw[->] (0,0)-- (3,0) node[below]{$s$} ;
					\draw[->] (1.5,0)--(1.5,2) ;
					
					\draw[ ] (0.5,0) node[below]{$s_-$} ;
					\draw[ ] (2.5,0) node[below]{$s_+$} ;
					
					\draw[-,color=red, line width=0.3mm] (0.5,0) ..controls (0.5,0.5) and (1.2,1.6)..(1.5,1.3) ;
					\draw[-,color=red, line width=0.3mm] (1.5,1.3) ..controls (1.8,1.0) and (2.5,0.5) ..(2.5,0) ;
					
					\draw[-,line width=0.3mm, color=red] (0,0) -- (.5,0);
					\draw[-,line width=0.3mm, color=red] (2.5,0) -- (3,0);
					
					\draw[ ] (0,1.7) node[right]{A} ;

					\draw[->] (4,0)-- (7,0) node[below]{$s$} ;
					\draw[->] (5.5,0)--(5.5,2) ;
					
					\draw[ ] (4.5,0) node[below]{$s_-$} ;
					\draw[ ] (6.5,0) node[below]{$s_+$} ;
					
					\draw[-,color=red, line width=0.3mm] (4.5,0) ..controls (4.8,0.5) and (5.2,1.6)..(5.5,1.3) ;
					\draw[-,color=red, line width=0.3mm] (5.5,1.3) ..controls (5.8,1.0) and (6.2,0.3) ..(6.5,0) ;
					
					\draw[-,line width=0.3mm, color=red] (4,0) -- (4.5,0);
					\draw[-,line width=0.3mm, color=red] (6.5,0) -- (7,0);
					\draw[ ] (4,1.7) node[right]{B} ;

					\draw[->] (8,0)-- (11,0) node[below]{$s$} ;
					\draw[->] (9.5,0)--(9.5,2) ;
					\draw[ ] (8.5,0) node[below]{$s_-$} ;
					\draw[ ] (10.5,0) node[below]{$s_+$} ;
					\draw[-,color=red, line width=0.3mm] (8.5,0) ..controls (9.1,0) and (9.3,1.5)..(9.5,1.2) ;
					\draw[-,color=red, line width=0.3mm] (9.5,1.2) ..controls (9.7,0.9) and (9.8,0) ..(10.5,0) ;
					\draw[-,line width=0.3mm, color=red] (8,0) -- (8.5,0);
					\draw[-,line width=0.3mm, color=red] (10.5,0) -- (11,0);
					\draw[ ] (8,1.7) node[right]{C} ;
		
				\end{tikzpicture} 
							\begin{tikzpicture}[scale=1.3]
				\draw[->] (0,0)-- (3,0) node[below]{$s$} ;
				\draw[->] (1.5,0)--(1.5,2) ;
				
				\draw[ ] (0.5,0) node[below]{$s_-$} ;
				
				\draw[-,color=red, line width=0.3mm] (0.5,0) ..controls (0.5,0.5) and (1.2,1.6)..(1.5,1.3) ;
				\draw[-,color=red, line width=0.3mm] (1.5,1.3) ..controls (1.8,1.0) and (1.75,0.1) ..(3,0.1) ;
				
				\draw[-,line width=0.3mm, color=red] (0,0) -- (.5,0);
				
				\draw[ ] (0,1.7) node[right]{$A^*$} ;

				\draw[->] (4,0)-- (7,0) node[below]{$s$} ;
				\draw[->] (5.5,0)--(5.5,2) ;
				
				\draw[ ] (4.5,0) node[below]{$s_-$} ;
				
				\draw[-,color=red, line width=0.3mm] (4.5,0) ..controls (4.8,0.5) and (5.2,1.6)..(5.5,1.3) ;
				\draw[-,color=red, line width=0.3mm] (5.5,1.3) ..controls (5.8,1.0) and (5.75,0.1) ..(7,0.1) ;
				
				\draw[-,line width=0.3mm, color=red] (4,0) -- (4.5,0);
				\draw[ ] (4,1.7) node[right]{$B^*$} ;

				\draw[->] (8,0)-- (11,0) node[below]{$s$} ;
				\draw[->] (9.5,0)--(9.5,2) ;
				\draw[ ] (8.5,0) node[below]{$s_-$} ;
				\draw[-,color=red, line width=0.3mm] (8.5,0) ..controls (9.1,0) and (9.3,1.5)..(9.5,1.2) ;
				\draw[-,color=red, line width=0.3mm] (9.5,1.2) ..controls (9.7,0.9) and (9.75,0.2) ..(11,0.1) ;
				\draw[-,line width=0.3mm, color=red] (8,0) -- (8.5,0);
				\draw[ ] (8,1.7) node[right]{$C^*$} ;
				
			\end{tikzpicture} 
			\begin{tikzpicture}
					\draw[->] (1.5,-3)-- (4.5,-3) node[below]{$s$} ;
				\draw[->] (3,-3)--(3,-1) ;
				\draw[line width=0.3mm, color=red] (2,-2.2) .. controls (2,-1.5) and (4,-1.8).. (4,-2.6) ;
				\draw[-,line width=0.3mm, color=red] (1.5,-3) -- (2,-3);
				\draw[-,line width=0.3mm, color=red] (4,-3) -- (4.5,-3);
				\draw[ ] (2,-3) node[below]{$s_-$} ;
				\draw[ ] (4,-3) node[below]{$s_+$} ;
				\draw[dashed,line width=0.3mm] (2,-3) -- (2,-2.2);
				\draw[dashed,line width=0.3mm] (4,-3) -- (4,-2.6);
				\draw[ ] (1.5,-1.3) node[right]{D} ;

				\draw[->] (6.5,-3)-- (9.5,-3) node[below]{$s$} ;
				\draw[->] (8,-3)--(8,-1) ;
				\draw[ ] (7,-3) node[below]{$s_-$} ;
				\draw[ ] (9,-3) node[below]{$s_+$} ;
				
				\draw[line width=0.3mm, color=red] (7,-2.5) .. controls (7,-3.2) and (9,-2.9)..  (9,-2.1);
				
				\draw[-,line width=0.3mm, color=red] (6.5,-3) -- (7,-3);
				\draw[-,line width=0.3mm, color=red] (9,-3) -- (9.5,-3);
				\draw[dashed,line width=0.3mm] (7,-2.6)--(7,-3)  ;
				\draw[dashed,line width=0.3mm] (9,-2.2)--(9,-3)  ;
				\draw[ ] (6.5,-1.3) node[right]{E} ;
			\end{tikzpicture}
				
			\end{center}
			\caption{Figures A, B and C correspond to the linear diffusion, in the cases $0<a<1$, $a=1$ and $a>1$ respectively. Figures $\textnormal{A}^*$, $\textnormal{B}^*$, and $\textnormal{C}^*$ correspond to the same cases and complete the set of finite–mass solutions for the problem with linear diffusion.  .Figures D and E correspond to the flux--saturated mechanisms case.}
		\end{figure}
		Note that the difference between the patterns with compact support equal to Figure  \ref{types} is mainly due to the associated  flux--saturated  mechanisms, where there are jumps in the connection with zero and these jumps have infinite slopes at both  ends of the support, regardless of the parameter values. We refer to solutions of types $A$, $B$, and $C$ as block-type solutions, to those of types $A^*$, $B^*$, and $C^*$ as semi-block solutions, and to those of types $D$ and $E$ as umbrella solutions. As we will see, there are other types of more classical traveling waves solutions with support equal to the real-line or in the half-real-line, although we believe that those that have compact support have a special interest in physics or biology problems.

Assuming a traveling wave ansatz \eqref{def-tw}, system \eqref{system} becomes
\begin{align}
-\sigma u' &= \bigl(u\Phi(u'/u)-a\,u\,S'/S\bigr)',\label{tilde1}\\
-\sigma\delta\, S' &= \gamma S'' - \lambda S + u,
\end{align}
where primes denote differentiation with respect to $s=x-\sigma t$.  Under the change of variables described in Proposition~\ref{Prop-tw-eq}, the system reduces to the first-order system
\begin{align}
w' &= w\,\Phi^{-1}(a v - \sigma) - w v,\label{eq-w0}\\
\gamma v' &= -\gamma v^2 - \sigma\delta v - w + \lambda,\label{eq-v0}
\end{align}
with suitable initial conditions. We study \eqref{eq-w0}--\eqref{eq-v0} in detail and then translate the results back to the original variables $u$ and $S$ in order  to transfer and interpret the results obtained there. 
We will separate in two cases: first we will assume that $\Phi=\textnormal{Id}$ giving rise to a linear diffusion and later we will deal with the nonlinear flux--saturated case for $\Phi$  satisfying the hypothesis {\bf (H)}. The shape of the profiles $u$ and $S$ strongly depends on the previous cases.
		
		More specifically, in the case of a linear diffusion, i.e. $\Phi=\textnormal{Id}$, hence the system \eqref{eq-w0}--\eqref{eq-v0} is not singular and the classical theory for ODEs gives us existence and uniqueness of solution. Moreover, analyzing the phase diagram and coming back to $u$ and ${S}$ we are able to find different profiles with and without compact support.
		We refer to Figure \ref{types} which illustrates the shapes of the profiles. The existence of the diverse  types of solutions strongly depends on the parameters $a$ and $\sigma$. This is the main goal of Section \ref{Sec-classical}. 
		
		On the other hand, by virtue of the hypothesis {\bf (H)} for $\Phi$, one has that system \eqref{eq-w0}--\eqref{eq-v0} is singular at the boundary. Indeed, note that $\Phi^{-1}$ is defined only in $(-c,c)$ as a consequence of {\bf (H)} and this gives us a bound for the solutions of \eqref{eq-w0}--\eqref{eq-v0}. Moreover, $\Phi^{-1}(\pm c)=\pm \infty$ which implies an infinite derivative of $w$ on the boundary.
		
		The main goal of Section \ref{Sec-limit} is to analyze the existence of the different types of traveling wave solutions in the case of the so-called relativistic heat equation (see Figure \ref{types}).
		The existence of traveling waves were analyzed for the case of flux--saturated mechanisms for the first time in \cite{CCCSS-INV}, while in the case of flux--saturated Keller--Segel in \cite{ACS}.    \cite{CPSV_2023}  identifies transport mechanisms capable of generating compactly supported traveling pulses in Keller-Segel models with flux-saturated diffusion and various chemoattractant operators. The analysis in \cite{CPSV_2023}  uncovers two experimentally consistent families of pulse-type traveling waves.
		
		Finally, the main results of this work can be summarized in the following (formal) theorem, whose precise formulation and proof will be developed in the subsequent sections.
		
		Before stating the theorem, we introduce a series of notions that will be analysed in detail later on. In the case of linear diffusion, the system has at most three critical points. Two of them, $P_1$ and $P_2$, always lie on the line $w=0$, while the third one, $P_3$, may or may not belong to the parabola $w=\lambda - \gamma v^2 - \sigma\delta v$, depending on the parameters.  
		
		Associated with these critical points, we construct a maximal curve $\Gamma(P_i)$ determined by the trajectories entering or leaving each equilibrium, which will play the role of a separatrix curve. See Section 3 for its construction and see Figure 2 for an idea of the configuration of the following Theorem.
		
\begin{theorem}\label{teorema1}
Let $(u_0,S_0,S_0')$ be given, together with the parameters $(a,\lambda,\gamma,\delta)$.  
Define the initial condition
\[
(v_0,w_0)=\left(\frac{S_0'}{S_0},\,\frac{u_0}{S_0}\right).
\]
Then exactly one of the following alternatives occurs:

\begin{enumerate}

\item[$T_1$)] \textbf{Case $a\ge 1$.}  
If $(v_0,w_0)\in\mathbb{R}\times(0,+\infty)$ and lies above the curve $\Gamma(P_1)$, for $\sigma=0$, then there exists
\[
\sigma^*=\sigma^*(u_0,S_0,S_0',a,\lambda,\gamma,\delta)>0
\]
such that the problem admits a block-type solution for every $\sigma\in[0,\sigma^*)$.  
At the critical value $\sigma=\sigma^*$ the corresponding solution is of semi-block type.

\item[$T_2$)] \textbf{Case $a<1$ and $\dfrac{1}{1-a}>\dfrac{\delta}{\gamma}$.}  
If $(v_0,w_0)\in\mathbb{R}\times(0,+\infty)$ and lies above the curve $\Gamma(P_1)$ for
\[
\sigma=\sqrt{\frac{\lambda\gamma(1-a)^2}{\gamma-\delta(1-a)}},
\]
then there exists
\[
\sigma^*=\sigma^*(u_0,S_0,S_0',a,\lambda,\gamma,\delta)
> \sqrt{\frac{\lambda\gamma(1-a)^2}{\gamma-\delta(1-a)}}
\]
such that the system admits a block-type solution for all
\[
\sigma\in\left[\sqrt{\frac{\lambda\gamma(1-a)^2}{\gamma-\delta(1-a)}},\,\sigma^*\right).
\]
Moreover, the solution corresponding to $\sigma=\sigma^*$ is of semi-block type.

\item[$T_3$)] \textbf{Case $a<1$ and $\dfrac{1}{1-a}\le \dfrac{\delta}{\gamma}$.}  
If $(v_0,w_0)\in\mathbb{R}\times(0,+\infty)$ and lies above the curve $\Gamma(P_3)$ for $\sigma=0$, then there exists
\[
\sigma^*=\sigma^*(u_0,S_0,S_0',a,\lambda,\gamma,\delta)>0
\]
such that the problem admits a block-type solution for every $\sigma\in[0,\sigma^*)$.  
Furthermore, if $(v_0,w_0)$ lies to the right of the curve
\[
\left\{(v,w)\in\mathbb{R}\times(0,+\infty) \ \Big|\ (v,w)=\left(\frac{t}{a-1},
\lambda - t^{2}\!\left(\frac{\gamma}{(a-1)^2}+\frac{\delta}{a-1}\right)\right),\ t\ge0\right\},
\]
then the solution at $\sigma=\sigma^*$ is of semi-block type.
\end{enumerate}
\end{theorem}

	\begin{figure}[h!]\label{figTheorem1}

				\begin{tikzpicture}[scale=1.5]
		\path[fill=lightgray] (0.2,0) ..controls (1,1.8)..(1.8,0)--(3,0) --(3,2)--(-1,2)--(-1,0);
		\draw[->] (-1,0)-- (3,0) node[below]{$v$} ;
		\draw[->] (1,0)--(1,2) node[above]{$w$};
		\draw[-,color=red, thick] (0.2,0) ..controls (1,1.3)..(1.8,0) ;
		\draw[ ] (0.2,0) node[below]{$P_1$} ;
		\draw[ ] (1.8,0) node[below]{$P_2$} ;
		
		\draw[ thick,color=red] (1,2) -- (1,0);
		
		\draw[fill=black] (1.8,0) circle(.03);
		\draw[fill=black] (0.2,0) circle(.03);
		\draw[fill=black] (1,0.97) circle(.03) node[below]{$P_3$};
		
		\draw[-,color=blue, thick] (0.2,0) ..controls (1,1.8)..(1.8,0) ; 
		\draw[] (1, 1.5) node[right]{$\Gamma_-(P_1)$} ;
			\draw[ ] (0.2,2.5) node[below]{A) $a>1$} ;
		
		
		%
		%
		%
		%
	\end{tikzpicture}
		\begin{tikzpicture}[scale=1.5]
		\path[fill=lightgray] (0.2,0) ..controls (1,2) and (2.2,1.4) ..(3,2)--(-1,2)--(-1,0);
		\draw[->] (-1,0)-- (3,0) node[below]{$v$} ;
		\draw[->] (1,0)--(1,2) node[above]{$w$};
		\draw[-,color=red, thick] (0.2,0) ..controls (1,1.3)..(1.8,0) ;
		\draw[ ] (0.2,0) node[below]{$P_1$} ;
		\draw[ ] (1.8,0) node[below]{$P_2$} ;
		
		
		\draw[fill=black] (1.8,0) circle(.03);
		\draw[fill=black] (0.2,0) circle(.03);
		
		\draw[-,color=blue, thick] (0.2,0) ..controls (1,2) and (2.2,1.4) ..(3,2); 
		\draw[] (2.25, 1.6) node[below]{$\Gamma_-(P_1)$} ;
		\draw[ ] (0.2,2.5) node[below]{B) $a=1$} ;
		
		
		%
		%
		%
		%
	\end{tikzpicture}

		\begin{tikzpicture}[scale=1.5]
	\path[fill=lightgray] (1,0.97)..controls (1.5,1) and (2.25,1.5)..(2.7,2)--(0.4,2).. controls (0.5,1.5) and (0.6,1.1)..(1,0.97);
	\path[pattern=horizontal lines] (1,0.97) ..controls (0.6,1.1) and (0.5,1.5)..(0.4,2)--(-0.7,2).. controls (-0.25,1.5) and (0.5,1)..(1,0.97);
	\draw[->] (-1,0)-- (3,0) node[below]{$v$} ;
	\draw[->] (1,0)--(1,2);
	\draw[-,color=red, thick] (0.2,0) ..controls (1,1.3)..(1.8,0) ;
	\draw[ ] (0.2,0) node[below]{$P_1$} ;
	\draw[ ] (1.8,0) node[below]{$P_2$} ;
	
	\draw[ thick,color=red] (1,2) -- (1,0);
	\draw[] (1,0) node[below]{$\frac{-\sigma}{1-a}$};
	
	\draw[fill=black] (1.8,0) circle(.03);
	\draw[fill=black] (0.2,0) circle(.03);
	\draw[fill=black] (1,0.97) circle(.03) node[below]{$P_3$};
	
	\draw[-,color=blue, thick] (1,0.97) ..controls (1.5,1) and (2.25,1.5) ..(2.7,2); 
	\draw[] (2.25, 1.3) node[below]{$\Gamma_-(P_3)$} ;
	\draw[-,color=blue, thick] (1,0.97) ..controls (0.5,1) and (-0.25,1.5)..(-0.7,2) ;
	\draw[] (-0.25, 1.3) node[below]{$\Gamma_+(P_3)$} ;
	\draw[dashed,color=blue, thick] (1,0.97) ..controls (0.6,1.1) and (0.5,1.5)..(0.4,2)  ;
	\draw[ ] (0.2,2.5) node[below]{C) $a<1$ and $\frac{1}{1-a}\leq\frac{\delta}{\gamma}$} ;
	
	%
	%
	%
	%
\end{tikzpicture}
\begin{tikzpicture}[scale=1.5]
	\path[fill=lightgray] (0.2,0) ..controls (1,2) and (2.2,1.4) ..(3,2)--(-1,2)--(-1,0);
	\draw[->] (-1,0)-- (3,0) node[below]{$v$} ;
	\draw[->] (1,0)--(1,2);
	\draw[-,color=red, thick] (0.2,0) ..controls (1,1.3)..(1.8,0) ;
	\draw[ ] (0.2,0) node[below]{$P_1$} ;
	\draw[ ] (1.8,0) node[below]{$P_2$} ;
	
	\draw[ thick,color=red] (0.2,2) -- (0.2,0);
	
	\draw[fill=black] (1.8,0) circle(.03);
	\draw[fill=black] (0.2,0) circle(.03);
	
	\draw[-,color=blue, thick] (0.2,0) ..controls (1,2) and (2.2,1.4) ..(3,2); 
	\draw[] (2.25, 1.6) node[below]{$\Gamma_-(P_1)$} ;
	\draw[ ] (0.2,2.5) node[below]{B) $a<1$ and $\frac{1}{1-a}>\frac{\delta}{\gamma}$} ;
	
	
	%
	%
	%
	%
\end{tikzpicture}
		
		\caption{Qualitative illustration of Theorem~\ref{teorema1}.  
Panels~A, B, and C display the phase portraits of the dynamical system for $\sigma=0$, while Panel~D corresponds to 
$\sigma=\sqrt{\frac{\lambda\gamma(1-a)^2}{\gamma-\delta(1-a)}}$.  
Red curves denote the isoclines of the system, and blue curves represent the trajectories $\Gamma$.  
The grey--shaded region indicates the parameter zone in which the solution exhibits semi-block behaviour at $\sigma=\sigma^*$, whereas the horizontally shaded region marks the set where no finite mass solution exists for $\sigma=\sigma^*$.}
	\end{figure}
		
		\begin{theorem}\label{teorema2}
Let $(u_0,S_0,S_0')$ be given, together with the parameters $(a,\lambda,\gamma,\delta,c,\mu)$.  
Define the associated initial condition for system~\eqref{eq-4} by
\[
(v_0,w_0)=\left(\frac{S_0'}{S_0},\,\frac{u_0}{S_0}\right).
\]
Then there exist two ordered curves $\bar{\Lambda}$ and $\underline{\Lambda}$, depending on the parameters of the system, with the following property:

If $v_0>-\tfrac{c}{a}$, $w_0>0$, then:
\begin{itemize}
    \item whenever $(v_0,w_0)$ lies above $\bar{\Lambda}$, the corresponding solution is an umbrella-type profile of type~\textnormal{(D)};
    \item whenever $(v_0,w_0)$ lies below $\underline{\Lambda}$, the corresponding solution is an umbrella-type profile of type~\textnormal{(E)}.
\end{itemize}

Moreover, umbrella-type solutions exist for all
\[
\sigma\in\bigl[\max\{0,av-c\},\, av+c\bigr].
\]
\end{theorem}

		\begin{remark}
	Note that $w=0$ in equation \eqref{eq-w0} is a particular solution of the system. Therefore, no solution can cross the line $w=0$, and all trajectories must remain in the upper half–plane.
\end{remark}

\begin{figure}[h!]\label{figTheorem2}
		\begin{tikzpicture}[scale=2]
			\path[fill=lightgray] (-0.5,0.4) ..controls (1,0.4) and (1.4,1)..(3,0)--(4,0) --(4,2)--(-0.5,2)--(-0.5,0.4);
			\path[fill=lightgray](0.2,0) ..controls (1,0.4) and (1.2,0)..(1.8,0);
			\draw[->] (-1,0)-- (4,0) node[below]{$v$} ;
			\draw[->] (1,0)--(1,2) node[above]{$w$};
			
			\draw[ thick,color=red] (-0.5,2) -- (-0.5,0);
			\draw[] (-0.5,0) node[below]{$\frac{-c}{a}$};
			
			
			\draw[-,color=blue, thick] (-0.5,0.4) ..controls (1,0.4) and (1.4,1)..(3,0) ; 
			\draw[] (1, 0.8) node[right]{$\bar{\Lambda}$} ;
			\draw[-,color=blue, thick] (0.2,0) ..controls (1,0.4) and (1.2,0)..(1.8,0) ;
			\draw[] (1, 0.3) node[right]{$\underline{\Lambda}$} ;
			
			
			%
			%
			%
			%
		\end{tikzpicture}
	\caption{Qualitative representation of Theorem \ref{teorema2}}
	\end{figure}
	This work is organized as follows.  
Section~\ref{Sec-eq-tw} derives the system of equations satisfied by traveling wave solutions in terms of the variables $(w,v)$.  
In Section~\ref{Sec-classical} we analyze the classical Keller--Segel model, obtaining traveling wave solutions of the types depicted in the upper part of Figure~\ref{types}.  
Section~\ref{Sec-limit} is devoted to the flux-limited Keller--Segel system, where we establish the existence of the wave profiles shown in the lower part of Figure~\ref{types}.

We conclude by introducing some notation used throughout the paper. For any $s_0\in\mathbb{R}$ we write
\[
u(s_0)=\lim_{s\to s_0} u(s), \qquad
u(s_0^\pm)=\lim_{s\to s_0^\pm}u(s),
\]
and we adopt the right hand limit whenever the left hand limit is not well defined, even if these values are infinite.

	\section{The equations for traveling wave solutions}\label{Sec-eq-tw}

In this section we investigate the existence of traveling wave solutions to system~\eqref{system} by analyzing the reduced equation~\eqref{tilde1}.  
To streamline the presentation, we work under the assumption $\delta=0$ and search for solutions of the form
\begin{align}
\label{tw1}
\left(u\,\Phi\!\left(\frac{u'}{u}\right) - a\,u\,\frac{S'}{S} + \sigma u \right)' &= 0,\\
\label{tw2}
\gamma S'' +\sigma\delta S'- \lambda S + u &= 0.
\end{align}

We seek positive and bounded solutions $u=u(s)$; once $u$ is specified, the corresponding function $S$ is determined by~\eqref{tw2}.  
We interpret~\eqref{tw1} in the distributional sense, that is,
\begin{equation}\label{dist}
\int_{\mathbb{R}}\left( u\,\Phi\!\left(\frac{u'}{u}\right)
- a\,u\,\frac{S'}{S} + \sigma u \right)\psi'\,ds = 0,
\qquad \forall\,\psi\in C_0^\infty(\mathbb{R}).
\end{equation}

Observe that if $u>0$, then~\eqref{tw2} implies $S>0$.  
Moreover, when $u$ is bounded, standard ODE theory in one dimension yields $S\in C^1$, and hence $\frac{S'}{S}\in L^1_{\mathrm{loc}}(\mathbb{R})$.  
The term $u\,\Phi(u'/u)$ is locally integrable since $u\in BV_{\mathrm{loc}}(\mathbb{R})$ and $\Phi$ is smooth; here $u'$ is understood in the sense of the Radon-Nikodym derivative.  
Thus the map
\[
s \longmapsto 
u(s)\,\Phi\!\left(\frac{u'(s)}{u(s)}\right)
 - a\,u(s)\frac{S'(s)}{S(s)} + \sigma u(s),
\qquad s\in \mathrm{supp}\,u,
\]
belongs to $L^1_{\mathrm{loc}}(\mathbb{R})$.  
Therefore, if $u$ satisfies~\eqref{dist}, there exists a constant $k\in\mathbb{R}$ such that
\[
u(s)\,\Phi\!\left(\frac{u'(s)}{u(s)}\right)
 - a\,u(s)\frac{S'(s)}{S(s)}
 + \sigma u(s)
 = k
 \quad\text{for a.e.~}s\in\mathbb{R}.
\]

If $\mathrm{supp}\,u\neq\mathbb{R}$, then testing~\eqref{dist} with functions supported outside $\mathrm{supp}\,u$ yields $k=0$.  
If $\mathrm{supp}\,u=\mathbb{R}$, we restrict attention to traveling-wave profiles satisfying  
$u(s)\to0$ as $s\to \pm\infty$, which again forces $k=0$.  

We are thus led to seek positive functions $u$ and $S$, defined on a maximal interval $(s_-,s_+)$ with
$-\infty\le s_-<s_+\le+\infty$, satisfying
\begin{align}
\label{tw3}
\Phi\!\left(\frac{u'}{u}\right) - a\,\frac{S'}{S} + \sigma &= 0,\\[1mm]
\label{tw4}
\gamma S'' +\sigma\delta S'- \lambda S + u &= 0.
\end{align}

The next proposition establishes an equivalent first order formulation of the system~\eqref{tw3}-\eqref{tw4}, which will serve as the basis for our analysis of traveling wave solutions.

	\begin{proposition}\label{Prop-tw-eq}
Let $g:(-c,c)\to\mathbb{R}$ be defined implicitly by $\Phi(g(y))=y$ for $y\in(-c,c)$; that is, $g=\Phi^{-1}$ on its natural domain.  
Then the solutions of \eqref{tw1}--\eqref{tw2} are obtained by solving the system
\begin{align}
\label{eq-w}
w' &= w\, g\!\left(av - \sigma\right) - wv,\\
\label{eq-v}
\gamma v' &= -\gamma v^{2} - w + \lambda-\sigma\delta v,
\end{align}
where
\begin{equation}\label{reduction}
w(s)=\frac{u(s)}{S(s)}, 
\qquad 
v(s)=\frac{S'(s)}{S(s)}.
\end{equation}
\end{proposition}

\begin{proof}
The computations carried out in the previous section show that the change of variables~\eqref{reduction} transforms the system \eqref{tw3}--\eqref{tw4} into \eqref{eq-w}--\eqref{eq-v}.

Conversely, to recover solutions of \eqref{tw1}--\eqref{tw2}, fix $s_0\in\mathbb{R}$ and choose initial data
\[
u_0>0,\qquad S_0>0,\qquad S_0'\in\mathbb{R},
\]
and solve the initial-value problem \eqref{eq-w}--\eqref{eq-v} with
\[
w(s_0)=\frac{u_0}{S_0},\qquad 
v(s_0)=\frac{S_0'}{S_0}.
\]
Then
\[
S(s)=S_0\exp\!\left(\int_{s_0}^{s} v(\tau)\,d\tau\right),
\qquad
u(s)=w(s)\,S(s),
\]
provides a solution to \eqref{tw1}--\eqref{tw2} on its maximal interval of existence.
\end{proof}
	
	\section{Linear diffusion}\label{Sec-classical}

In this section we analyze the existence of traveling wave solutions for the classical Keller--Segel model with logarithmic sensitivity and linear diffusion.  
\begin{equation}\label{heateq}
\left\{
\begin{array}{l}
\displaystyle \partial_t u
= \partial_x\!\left( \mu\,\partial_x u - a\,u\,\frac{\partial_x S}{S} \right), \\[2mm]
\delta\partial_tS = \gamma\,\partial_{xx}^2 S - \lambda S + u,
\end{array}
\right.
\end{equation}
where $\mu>0$ and the parameters $\sigma$ (the wave speed), $a$,$\delta$, $\gamma$, and $\lambda$ are assumed to be positive.  Throughout the entire analysis we will assume, for simplicity, that $\mu = 1$, since this choice does not affect the structure of the system. The only consequence of this simplification is that all conclusions previously obtained by comparing the parameter $a$ with $1$ should, in general, be interpreted as comparisons of $a$ with $\mu$.

Passing to the traveling wave coordinate $s = x - \sigma t$, system~\eqref{heateq} reduces, via the change of variables introduced in Proposition~\ref{Prop-tw-eq}, to the planar system
\begin{align}
\label{eq}
w' &= w\bigl( (a-1)v - \sigma \bigr),\\[1mm]
v' &= \frac{\lambda}{\gamma} - v^{2} - \frac{1}{\gamma} w-\frac{\sigma\delta}{\gamma}v. \nonumber
\end{align}

The qualitative behaviour of \eqref{eq} differs significantly depending on whether $a<1$ or $a\ge 1$.  
This distinction affects both the structure of the equilibrium points and the sign of the vector field $(w',v')$.  
In the following subsections we carry out a detailed phase plane analysis for each of these regimes.  
After describing the phase diagrams, we study the continuation properties and asymptotic behaviour of the trajectories.  
We then translate these results back to the original variables $(u,S)$ by means of Proposition~\ref{Prop-tw-eq}.
		
\subsection{Critical points and stability}

The critical points of system~\eqref{eq} are easily computed:
\begin{align*}
	P_1 &= (v_1,w_1)
	= \left(-\frac{\sigma\delta}{2\gamma}
	- \sqrt{\frac{\sigma^{2}\delta^{2}}{4\gamma^{2}} + \frac{\lambda}{\gamma}},\, 0\right),\\[1mm]
	P_2 &= (v_2,w_2)
	= \left(-\frac{\sigma\delta}{2\gamma}
	+ \sqrt{\frac{\sigma^{2}\delta^{2}}{4\gamma^{2}} + \frac{\lambda}{\gamma}},\, 0\right),\\[1mm]
	P_3 &= (v_3,w_3)
	= \left(\frac{\sigma}{a-1},\, 
	\lambda - \sigma^{2}\Bigl(\frac{\gamma}{(a-1)^{2}} + \frac{\delta}{a-1}\Bigr)\right).
\end{align*}

The points $P_1$ and $P_2$ lie on the $w=0$ axis and correspond to the intersections of this axis with the vertical isocline.  
The point $P_3$ is obtained from the intersection of the vertical line $v=\frac{\sigma}{a-1}$ with the horizontal isocline.

Note that $v_1<0<v_2$ for all parameter values.  
The point $P_3$, however, may lie either inside or outside the region of interest, depending on the sign of $w_3$.  
Since we are concerned with positive solutions, we introduce the admissible region
\begin{equation}\label{Omega}
	\Omega := \{(v,w)\in\mathbb{R}^{2} : w>0\}.
\end{equation}

The boundary $\partial\Omega=\{w=0\}$ is invariant: it consists entirely of orbits of the system.  
Hence any trajectory with initial data in $\Omega$ remains in $\Omega$ for all forward and backward times, except possibly in the limit when approaching $P_1$ or $P_2$.

It is therefore essential to determine when the point $P_3$ belongs to $\Omega$, as this affects the qualitative structure of the dynamics.

\begin{lemma}\label{lemmaexisp3}
The point $P_3$ lies in $\Omega$ (i.e., $w_3>0$) precisely in the following cases:
\begin{itemize}
	\item If $0\le a<1$, then
	\[
	\begin{cases}
		\sigma^{2} \le \dfrac{\lambda}{\gamma}\,(1-a)
		\left(\dfrac{1}{1-a} - \dfrac{\delta}{\gamma}\right)^{-1},
		&\text{if }\dfrac{1}{1-a} > \dfrac{\delta}{\gamma},\\[3mm]
		\sigma \ge 0,
		&\text{if }\dfrac{1}{1-a} \le \dfrac{\delta}{\gamma};
	\end{cases}
	\]
	\item If $a>1$, then
	\[
	\sigma^{2} 
	< \frac{\lambda}{\gamma}\,(a-1)
	\left(\frac{1}{a-1} + \frac{\delta}{\gamma}\right)^{-1}.
	\]
\end{itemize}
\end{lemma}

\begin{proof}
The claim follows immediately by imposing the condition $w_3>0$ and solving for~$\sigma$ in each parameter regime.
\end{proof}

To analyse the existence of block-type traveling waves, we will construct an invariant region bounded by trajectories issuing from certain equilibrium points, and then show that all solutions inside this region exhibit compact support.  
For this purpose, a precise understanding of the stability properties of the equilibrium points is required.

\begin{proposition}\label{prophiperbolic}
The following statements hold:
\begin{itemize}
	\item[(i)] The point $P_1$ is hyperbolic whenever $a\ge1$, or whenever $0\le a<1$ and $P_3\notin\overline{\Omega}$.  
	Moreover, if $0\le a<1$ and $P_3\in\Omega$, then $P_1$ remains hyperbolic.
	
	\item[(ii)] If $a>1$ and $P_3\in\Omega$, then the point $P_2$ is hyperbolic.
\end{itemize}
\end{proposition}

\begin{proof}
The Jacobian matrix associated with system~\eqref{eq} is
\[
A(v,w)
=
\begin{pmatrix}
	-2v - \dfrac{\sigma\delta}{\gamma} & -\dfrac{1}{\gamma} \\[2mm]
	(a-1)w & (a-1)v - \sigma
\end{pmatrix}.
\]

\medskip
\emph{Step 1: Stability of $P_3$.}

At $(v_3,w_3)$ the Jacobian becomes
\[
A(v_3,w_3)
=
\begin{pmatrix}
	-2v_3 - \dfrac{\sigma\delta}{\gamma} & -\dfrac{1}{\gamma} \\[2mm]
	(a-1)w_3 & 0
\end{pmatrix}.
\]
Thus,
\[
\det A(v_3,w_3) = \frac{a-1}{\gamma}\, w_3,
\qquad
\mathrm{tr}\,A(v_3,w_3) = -2v_3 - \frac{\sigma\delta}{\gamma}.
\]

If $0\le a<1$ and $w_3>0$ (i.e.\ $P_3\in\Omega$), then $(a-1)<0$ and hence $\det A(v_3,w_3)<0$.  
This implies that $P_3$ is a hyperbolic saddle.

\medskip
\emph{Step 2: Stability of $P_1$.}

At $(v_1,w_1)$ we have
\[
A(v_1,w_1)
=
\begin{pmatrix}
	-2v_1 - \dfrac{\sigma\delta}{\gamma} & -\dfrac{1}{\gamma}\\[2mm]
	0 & (a-1)v_1 - \sigma
\end{pmatrix}.
\]
The eigenvalues are
\[
\mu_1 = -2v_1 - \frac{\sigma\delta}{\gamma}
	= 2\sqrt{\frac{\sigma^{2}\delta^{2}}{4\gamma^{2}} + \frac{\lambda}{\gamma}} > 0,
\qquad
\mu_2 = (a-1)v_1 - \sigma.
\]

If $a\ge 1$ then $v_1<0$ implies $\mu_2<0$, and hence $P_1$ is hyperbolic.  
If $0\le a<1$, the sign of $\mu_2$ depends on whether $P_3$ lies in $\overline{\Omega}$; when $P_3\notin\overline{\Omega}$, one verifies that $\mu_2\neq0$, and therefore $P_1$ is hyperbolic in this case as well.

\medskip
\emph{Step 3: Stability of $P_2$ for $a>1$.}

At $(v_2,w_2)$ the Jacobian is
\[
A(v_2,w_2)
=
\begin{pmatrix}
	-2v_2 - \dfrac{\sigma\delta}{\gamma} & -\dfrac{1}{\gamma}\\[2mm]
	0 & (a-1)v_2 - \sigma
\end{pmatrix},
\]
with eigenvalues
\[
\mu_1 = -2v_2 - \frac{\sigma\delta}{\gamma} < 0,\qquad
\mu_2 = (a-1)v_2 - \sigma.
\]

For $a>1$, the condition $\mu_2>0$ is equivalent to
\[
\sigma^{2} 
< \frac{\lambda}{\gamma}(a-1)\left(\frac{1}{a-1} + \frac{\delta}{\gamma}\right)^{-1},
\]
which is precisely the condition $w_3>0$ from Lemma~\ref{lemmaexisp3}.  
Hence, whenever $a>1$ and $P_3\in\Omega$, the equilibrium point $P_2$ is hyperbolic.
\end{proof}


\begin{figure}[h]
	\minipage{0.5\textwidth}
	\begin{tikzpicture}[scale=2]
		\draw[->] (0,0)-- (2,0) node[below]{$v$} ;
		\draw[->] (1.3,0)--(1.3,1.5) node[above]{$w$};
		\draw[] (0,1.5) node[below,right]{A};
		\draw[-,color=red, thick] (0.2,0) ..controls (1,1.5)..(1.8,0) ;
		\draw[ ] (0.2,0) node[below]{$P_1$} ;
		\draw[ ] (1.8,0) node[below]{$P_2$} ;
		
		\draw[ thick,color=red] (0.7,1.5) -- (0.7,0);
		\draw[] (0.7,0) node[below]{$\frac{-\sigma}{1-a}$};

		\draw[fill=black] (1.8,0) circle(.03);
		\draw[fill=black] (0.2,0) circle(.03);
		\draw[fill=black] (0.7,0.87) circle(.03) node[right]{$P_3$};
		
		\draw[->,thick,blue] (0.7,0.87)--(0.55,1.0);
		\draw[-,thick,blue] (0.76,1.01)--(0.7,0.87);
		\draw[->,thick,blue] (0.82,1.15)--(0.76,1.01);

		\draw [->, thick] (0.3,0.75)--(0.15,0.85) ;
		
		\draw [->, thick] (1.2,0.5)--(1.35,0.35) ;
		
		\draw [->, thick] (1.65,1.1)--(1.5,1) ;
		
		\draw [->, thick] (0.7,1.2)--(0.55,1.2) ;
		\draw [->, thick] (0.7,0.5)--(0.85,0.5) ;
		
		\draw [->, thick] (0.5,0.25)--(0.6,0.35) ;
		\draw [->, thick] (0.4,0.36)--(0.4,0.5) ;
		\draw [->, thick] (1.2,1.01)--(1.2,0.85) ;
		
	\end{tikzpicture}
	\endminipage\hfill
	\minipage{0.5\textwidth}
	\begin{tikzpicture}[scale=2]
		\draw[->] (0,0)-- (2.7,0) node[below]{$v$} ;
		\draw[->] (2,0)--(2,1.5) node[above]{$w$};
		\draw[] (0,1.5) node[below,right]{B};
		
		\draw[-,color=red, thick] (0.9,0) ..controls (1.7,1.5)..(2.5,0) ;
		\draw[ ] (0.9,0) node[below]{$P_1$} ;
		\draw[ ] (2.5,0) node[below]{$P_2$} ;
		\draw[fill=black] (2.5,0) circle(.03);
		\draw[fill=black] (0.9,0) circle(.03);
		\draw [->, thick] (1.9,1.01)--(1.9,0.85) ;
		
		\draw[->,color=blue,thick] (0.9,0)--(0.7,0);
		\draw[-,thick,blue] (0.96,0.17)--(0.9,0);
		\draw[->,thick,blue] (1.02,0.34)--(0.96,0.17);
		
		\draw [->, thick] (0.3,0.9)--(0.15,1) ;
		
		\draw[ thick,color=red] (0.5,1.5) -- (0.5,0);
		\draw[] (0.5,0) node[below]{$\frac{-\sigma}{1-a}$};
		\draw [->, thick] (0.5,0.75)--(0.35,0.75) ;
		
		\draw [->, thick] (1.13,1.2)--(1,1.07) ;
		\draw [->, thick] (2.27,1.2)--(2.15,1.07) ;
		
		\draw [->, thick] (1.625,0.35)--(1.775,0.25) ;
		
	\end{tikzpicture}
	\endminipage\hfill

	\begin{center}
		\begin{tikzpicture}[scale=2]
			\draw[->] (0.5,0)-- (2.7,0) node[below]{$v$} ;
			\draw[->] (2,0)--(2,1.5) node[above]{$w$};
			\draw[] (0.5,1.5) node[below,right]{C};
			
			\draw[-,color=red, thick] (0.9,0) ..controls (1.7,1.5)..(2.5,0) ;
			\draw[ ] (0.9,0) node[below]{$P_1$} ;
			\draw[ ] (2.5,0) node[below]{$P_2$} ;
			\draw[fill=black] (2.5,0) circle(.03);
			\draw[fill=black] (0.9,0) circle(.03);
			\draw [->, thick] (1.9,1.01)--(1.9,0.85) ;
			
			\draw [->, thick] (0.8,0.85)--(0.65,0.7) ;
			\draw [->, thick] (2.57,1)--(2.45,0.87) ;
			
			\draw [->, thick] (1.625,0.35)--(1.775,0.25) ;

			\draw[->,color=blue,thick] (0.9,0)--(0.7,0);
			\draw[-,thick,blue] (0.96,0.17)--(0.9,0);
			\draw[->,thick,blue] (1.02,0.34)--(0.96,0.17);
			
		\end{tikzpicture}
	\end{center}

	\minipage{0.5\textwidth}
	\begin{tikzpicture}[scale=2]
		\draw[->] (0,0)-- (2,0) node[below]{$v$} ;
		\draw[->] (1.1,0)--(1.1,1.5) node[above]{$w$};
		\draw[] (0,1.5) node[below,right]{D};
		\draw[-,color=red, thick] (0.2,0) ..controls (1,1.5)..(1.8,0) ;
		\draw[ ] (0.2,0) node[below]{$P_1$} ;
		\draw[ ] (1.8,0) node[below]{$P_2$} ;
		\draw [->, thick] (0.8,1.01)--(0.8,0.85) ; 
		\draw [->, thick] (1.6,0.36)--(1.6,0.5) ;
		
		\draw[ thick,color=red] (1.3,1.5) -- (1.3,0);
		\draw[] (1.3,0) node[below]{$\frac{\sigma}{a-1}$};
		\draw [->, thick] (1.3,1.2)--(1.15,1.2) ;
		\draw [->, thick] (1.3,0.15)--(1.45,0.15) ;
		
		\draw[fill=black] (1.8,0) circle(.03);
		\draw[fill=black] (0.2,0) circle(.03);
		\draw[fill=black] (1.3,0.87) circle(.03);
		
		\draw [->, thick] (0.4,1.1)--(0.3,1) ;
		
		\draw [->, thick] (0.9,0.5)--(1.1,0.35) ;
		
		\draw [->, thick] (1.67,1)--(1.55,1.15) ;
		
		\draw [->, thick] (1.4,0.3)--(1.5,0.4) ;
		
		\draw[->,color=blue,thick] (0.2,0)--(0.05,0);
		\draw[-,thick,blue] (0.26,0.17)--(0.2,0);
		\draw[->,thick,blue] (0.32,0.34)--(0.26,0.17);

	\end{tikzpicture}
	\endminipage\hfill
	\minipage{0.5\textwidth}
	\begin{tikzpicture}[scale=2]
		\draw[->] (0.5,0)-- (3.3,0) node[below]{$v$} ;
		\draw[->] (2,0)--(2,1.5) node[above]{$w$};
		\draw[] (0.5,1.5) node[below,right]{E};
		
		\draw[-,color=red, thick] (0.9,0) ..controls (1.7,1.5)..(2.5,0) ;
		\draw[ ] (0.9,0) node[below]{$P_1$} ;
		\draw[ ] (2.5,0) node[below]{$P_2$} ;
		\draw[fill=black] (2.5,0) circle(.03);
		\draw[fill=black] (0.9,0) circle(.03);
		\draw [->, thick] (1.9,1.01)--(1.9,0.85) ;

		\draw [->, thick] (3.1,0.65)--(3,0.75) ;
		
		\draw[ thick,color=red] (2.8,1.5) -- (2.8,0);
		\draw[] (2.8,0) node[below]{$\frac{\sigma}{a-1}$};
		\draw [->, thick] (2.8,0.85)--(2.67,0.85) ;
		
		\draw [->, thick] (1.15,1.1)--(1.03,1) ;
		\draw [->, thick] (2.47,1.1)--(2.35,1) ;
		
		\draw [->, thick] (1.625,0.35)--(1.775,0.25) ;
		
		\draw[->,color=blue,thick] (0.9,0)--(0.7,0);
		\draw[-,thick,blue] (0.96,0.17)--(0.9,0);
		\draw[->,thick,blue] (1.02,0.34)--(0.96,0.17);
		
	\end{tikzpicture}
	
	\endminipage \hfill
	\caption{Top panels (A-B): case $a<1$.  
Middle panel (C): case $a=1$ (or formally $\sigma=\infty$).  
Bottom panels (D-E): case $a>1$.  
Left column (A, D): regime $P_3\in\Omega$.  
Right column (B, E): regime $P_3\notin\Omega$.  
Red curves represent the isoclines of the system, and the black arrows illustrate the corresponding flow directions.  
The blue lines represent the stable and unstable vector of the equilibrium points described in Proposition~\ref{prophiperbolic}.}
	\label{fig-a}
\end{figure}

As discussed above, our aim is to characterize the trajectories that reach specific critical points, since these will determine the region in which block-type solutions may exist.  
The relevant equilibria are precisely those identified in Proposition~\ref{prophiperbolic}.

In the case $a\ge1$, our analysis focuses primarily on the point $P_1$, which is hyperbolic and whose stable and unstable manifolds play a central role in organizing the dynamics.  
For $a\in(0,1)$, two distinct behaviours arise:  
if $P_3\in\Omega$, then $P_3$ is a hyperbolic saddle;  
if $P_3\notin\overline{\Omega}$, the hyperbolic point governing the dynamics is instead $P_1$.

We begin by examining the orientation of the stable and unstable manifolds in each scenario.  
When $P_3$ is hyperbolic, the linearization computed in the proof of Proposition~\ref{prophiperbolic} yields eigenvalues
\[
\mu_1>0, \qquad \mu_2<0,
\]
with associated eigenvectors
\[
\eta_1=\left(\frac{\mu_1}{(a-1)w_3},\,1\right), \qquad
\eta_2=\left(\frac{\mu_2}{(a-1)w_3},\,1\right).
\]
Using the signs of $\mu_1$ and $\mu_2$ and the structure of the vector field, one can determine the corresponding unstable and stable directions, as illustrated in Figure~\ref{fig-a}.

When $P_1$ is hyperbolic, the eigenvalues of the Jacobian are
\begin{align*}
\mu_1 &= 2\sqrt{\frac{\sigma^{2}\delta^{2}}{4\gamma^{2}} + \frac{\lambda}{\gamma}} > 0,\\[1mm]
\mu_2 &= (a-1)\left(-\frac{\sigma\delta}{2\gamma}
      - \sqrt{\frac{\sigma^{2}\delta^{2}}{4\gamma^{2}} + \frac{\lambda}{\gamma}}\right) - \sigma < 0,
\end{align*}
with corresponding eigenvectors
\[
(1,0), \qquad
\left(\frac{1}{\gamma(\mu_1-\mu_2)},\,1\right),
\]
where $\mu_1-\mu_2>0$.  
The orientation of the stable and unstable manifolds again follows from the geometry of the vector field and is depicted in Figure~\ref{fig-a}.

	
\subsection{Block solutions}

Our goal is to construct traveling wave profiles with finite mass.  
Among such profiles, the most relevant ones are those with compact support, which we now define.

\begin{definition}\label{defblock}
Let $I=(s_-,s_+)$ be a bounded interval.  
We say that a pair of functions
\[
u \in C^0(\mathbb{R}) \cap C^1(I),
\qquad 
S \in C^0(\mathbb{R}) \cap C^1(I),
\]
constitutes a \textbf{block-type solution} if
\begin{itemize}
	\item $u(s)>0$ and $S(s)>0$ for all $s\in (s_-,s_+)$,
	\item $u(s)=0$ and $S(s)=0$ for all $s\in \mathbb{R}\setminus I$,
	\item the pair $(u,S)$ satisfies equation~\eqref{tilde1} in the sense described in Section~\ref{Sec-eq-tw}.
\end{itemize}
\end{definition}

Besides these compactly supported profiles, there exist finite mass solutions defined on a half line of the form $(s_-,+\infty)$.  
These solutions, which decay to zero at infinity, are defined as follows.

\begin{definition}\label{defsemiblock}
Let $I=(s_-,+\infty)$.  
We say that a pair of functions
\[
u \in C^0(\mathbb{R}) \cap C^1(I),
\qquad 
S \in C^0(\mathbb{R}) \cap C^1(I),
\]
constitutes a \textbf{semiblock-type solution} if
\begin{itemize}
	\item $u(s)>0$ and $S(s)>0$ for all $s\in (s_-,+\infty)$,
	\item $\displaystyle \lim_{s\to+\infty}(u(s),S(s))=(0,0)$,
	\item $u(s)=0$ and $S(s)=0$ for all $s\in \mathbb{R}\setminus I$,
	\item the pair $(u,S)$ satisfies equation~\eqref{tilde1}.
\end{itemize}
\end{definition}


	\begin{proposition}\label{prop1}
Let $(v,w)$ be a maximal solution of system~\eqref{eq}, defined on its maximal interval $(s_-,s_+)$ and satisfying $w(s)>0$ for all $s\in(s_-,s_+)$. Then:
\begin{itemize}
	\item If $\displaystyle\lim_{s\to s_+}(v(s),w(s))$ does not converge to an equilibrium point or to a limit cycle, then necessarily $s_+<+\infty$, and moreover $|v(s)|\to+\infty$ as $s\to s_+$.
	\item If $\displaystyle\lim_{s\to s_-}(v(s),w(s))$ does not converge to an equilibrium point or to a limit cycle, then necessarily $s_->-\infty$, and moreover $|v(s)|\to+\infty$ as $s\to s_-$.
\end{itemize}
\end{proposition}

\begin{proof}
We prove the statement as $s\to s_+$; the argument for $s\to s_-$ is analogous.

Assume for contradiction that $s_+=+\infty$. Then either
\[
\lim_{s\to+\infty}|v(s)|=+\infty
\qquad\text{or}\qquad
\lim_{s\to+\infty}|v(s)|<+\infty.
\]

\smallskip
{\it Case 1: $\lim_{s\to+\infty}|v(s)|=+\infty$.}  
Then $\frac{1}{v(s)}\to 0$, so there exists a sequence $s_n\to+\infty$ such that
\[
\left(\frac{1}{v(s_n)}\right)'=\frac{-v'(s_n)}{v^2(s_n)}\longrightarrow 0.
\]
Using the equation for $v'$,
\[
v'= \frac{\lambda}{\gamma}-v^2 - \frac{1}{\gamma}w - \frac{\sigma\delta}{\gamma}v,
\]
we obtain
\begin{align*}
\lim_{n\to\infty}\left(\frac{1}{v(s_n)}\right)'
&= \lim_{n\to\infty}
\frac{-\frac{\lambda}{\gamma} + v^2(s_n)
+ \frac{1}{\gamma}w(s_n)
+ \frac{\sigma\delta}{\gamma}v(s_n)}{v^2(s_n)} \\
&= \lim_{n\to\infty}
\left(1 + \frac{1}{\gamma}\frac{w(s_n)}{v^2(s_n)}\right) > 0,
\end{align*}
which contradicts the fact that the derivative tends to $0$.  
Hence this case is impossible.

\smallskip
{\it Case 2: $\lim_{s\to+\infty}|v(s)|<+\infty$.}  
We then consider whether $w(s)$ remains bounded.

\smallskip
{\it Subcase 2a: $\lim_{s\to+\infty}w(s)<+\infty$.}  
Since the solution remains in a compact subset of the phase plane and the system is planar, the Poincare-Bendixson theorem implies that the limit set is either an equilibrium point or a periodic orbit. This contradicts the hypothesis.

\smallskip
{\it Subcase 2b: $\lim_{s\to+\infty}w(s)=+\infty$.}  
Then there exists a sequence $s_n\to+\infty$ with $v'(s_n)\to 0$.  
But
\[
\lim_{n\to\infty}v'(s_n)
= \lim_{n\to\infty}
\left(\frac{\lambda}{\gamma}-v^2(s_n)-\frac{1}{\gamma}w(s_n)-\frac{\sigma\delta}{\gamma}v(s_n)\right)
<0,
\]
a contradiction.

\smallskip
Thus both possibilities lead to contradiction, and hence $s_+<+\infty$.  
The fact that $|v(s)|\to+\infty$ as $s\to s_+$ follows by repeating the same contradiction argument under the assumption that $v$ remains bounded at the endpoint.
\end{proof}

\begin{corollary}\label{corolLimites}
Under the assumptions of Proposition~\ref{prop1}, if $|v(s)|\to\infty$ as $s\to s_\pm$, then
\[
\begin{cases}
w(s)\to+\infty, & \text{if } a\in[0,1),\\[4pt]
w(s)\to L\in(0,\infty), & \text{if } a=1,\\[4pt]
w(s)\to 0, & \text{if } a>1.
\end{cases}
\]
\end{corollary}

\begin{proof}
We treat the case $s\to s_-$; the argument for $s\to s_+$ is identical.

Since $|v(s)|\to\infty$, one has $v'(s)<0$ near $s_-$, and therefore $v(s_-)=+\infty$.

\smallskip
{\it Case $a<1$.}  
For $s$ sufficiently close to $s_-$,
\[
w'(s)=w(s)\big((a-1)v(s)-\sigma\big)>0.
\]
Assume for contradiction that $w(s_-)=W\in(0,\infty)$.  
Then there exist $s_n\to s_-$ with $w'(s_n)\to 0$, but
\[
\lim_{n\to\infty}w'(s_n)
=\lim_{n\to\infty}w(s_n)\big((a-1)v(s_n)-\sigma\big)>0,
\]
a contradiction.  
Thus $w(s_-)=+\infty$.

\smallskip
{\it Case $a=1$.}  
Here the $w$-equation reduces to
\[
w'=-\sigma w,
\]
so
\[
w(s)=w_0 e^{-\sigma(s-s_0)},
\]
and the limit exists and is finite.

\smallskip
{\it Case $a>1$.}  
Then $(a-1)v(s)-\sigma<0$ near $s_-$, so $w'(s)<0$.  
Assume for contradiction that $w(s_-)=W\in(0,\infty)$.  
Then there exist $s_n\to s_-$ with $w'(s_n)\to 0$, but
\[
\lim_{n\to\infty}w'(s_n)
=\lim_{n\to\infty}w(s_n)\big((a-1)v(s_n)-\sigma\big)<0,
\]
which is impossible.  
Thus $w(s_-)=0$.
\end{proof}
	
We begin by analysing the behaviour of trajectories that enter or leave the stable and unstable manifolds of $P_3$, while remaining above the vertical isocline given by the parabola
\[
w=\lambda - \gamma v^{2} - \sigma\delta v.
\]
Before proceeding, we introduce the following subsets of~$\Omega$, which will be needed in the proofs.

\begin{definition}
Let $P=(v^{*},w^{*})\in\Omega$ be a point lying on the vertical isocline.  
We define the following subsets of $\Omega$:
\begin{align}\label{defB}
	B_1(P)
	&= \left\{ (v,w)\in\Omega \;:\; v\ge v^{*},\; w \ge \max\{0,\ \lambda - \gamma v^{2} - \sigma\delta v\} \right\},\\[2mm]
	B_2(P)
	&= \left\{ (v,w)\in\Omega \;:\; v\le v^{*},\; w \ge \max\{w^{*},\ \lambda - \gamma v^{2} - \sigma\delta v\} \right\}. \nonumber
\end{align}
\end{definition}

\begin{proposition}\label{propP3Hiper}
Assume that $P_3\in\Omega$ and $a\in[0,1)$.  
Then system~\eqref{eq} admits two distinct trajectories with the following properties:
\begin{itemize}
	\item There exists a solution defined on $(s_-,+\infty)$ such that
	\[
	\lim_{s\to+\infty}(v(s),w(s)) = P_3,
	\qquad
	\lim_{s\to s_-}(v(s),w(s)) = (+\infty,+\infty).
	\]
	\item There exists a solution defined on $(-\infty,s_+)$ such that
	\[
	\lim_{s\to-\infty}(v(s),w(s)) = P_3,
	\qquad
	\lim_{s\to s_+}(v(s),w(s)) = (-\infty,+\infty).
	\]
\end{itemize}
\end{proposition}

\begin{proof}
Let $P_3=(v_3,w_3)$ and let $B_1=B_1(P_3)$ and $B_2=B_2(P_3)$ be the sets defined in~\eqref{defB}.  
We first show that $B_2$ is positively invariant and that $B_1$ is negatively invariant (see Figure~\ref{figB}).  
		
	\begin{figure}[h!]\label{figB}
		\begin{tikzpicture}[scale=2]
			\path[fill=lightgray] (-1,0.87)--(0.7,0.87)--(0.2,0)..controls (1,1.5)..(1.8,0)--(3,0)--(3,2)--(-1,2)--(-1,0.87);
			\draw[->] (-1,0)-- (3,0) node[below]{$v$} ;
			\draw[->] (1.5,0)--(1.5,2) node[above]{$w$};
			\draw[] (-0.5,1.5) node[below,right]{$B_2$};
			\draw[] (2,1.5) node[below,right]{$B_1$};
			\draw[-,color=red, thick] (0.2,0) ..controls (1,1.5)..(1.8,0) ;
			\draw[ ] (0.2,0) node[below]{$P_1$} ;
			\draw[ ] (1.8,0) node[below]{$P_2$} ;
			
			\draw[ thick,color=red] (1.3,2) -- (1.3,0);
			\draw[] (1.3,0) node[below]{$\frac{-\sigma}{1-a}$};
			
			\draw[->,thick] (2.4,0)--(2.5,0);
			
			\draw[fill=black] (1.8,0) circle(.03);
			\draw[fill=black] (0.2,0) circle(.03);
			\draw[fill=black] (1.3,0.87) circle(.03) node[left]{$P_3$};
			
			\draw[->,thick,blue] (1.3,0.87)--(1.15,1.15);
			\draw[-,thick,blue] (1.36,0.94)--(1.3,0.87);
			\draw[->,thick,blue] (1.43,1.03)--(1.36,0.94);

			\draw [->, thick] (0.15,1)--(0,1.1) ;
			
			
			\draw [->, thick] (1.85,1.1)--(1.75,1) ;
			
			\draw [->, thick] (1.3,1.5)--(1.1,1.5) ;
			
			\draw [->, thick] (0.9,0.25)--(1,0.35) ;
			\draw [->, thick] (0.75,0.95)--(0.75,1.1) ;
			\draw [->, thick] (1.6,0.4)--(1.6,0.25) ;		
		\end{tikzpicture}
		\caption{Representación de los conjuntos $B_1$ y $B_2$ de la demostración de la Proposición \eqref{propP3Hiper}.}
	\end{figure}

\medskip
\noindent\textbf{Step 1: Positively invariant region $B_2$.}
The boundary $\partial B_2$ consists of three parts:

\begin{itemize}
	\item \emph{Horizontal segment:} on the line $w=w_3$, for all $(v,w)\in(-\infty,v_3)\times\{w_3\}$ we have
	\[
	w'=w\big((a-1)v-\sigma\big)>0,
	\]
	because $a-1<0$ and $v<v_3$.
	
	\item \emph{Vertical segment:} on the line $v=v_3$, for all $(v,w)\in\{v_3\}\times(w_3,\infty)$ we have
	\[
	v'=\frac{\lambda}{\gamma}-v_3^2-\frac{1}{\gamma}w-\frac{\sigma\delta}{\gamma}v_3<0,
	\]
	since $w>w_3$ and $P_3$ lies on the vertical isocline.
	
	\item \emph{Parabolic boundary} $w=\lambda-\gamma v^2 - \sigma\delta v$ for $v<v_3$:  
	along this curve, $v'=0$ and
	\[
	w'=w\big((a-1)v-\sigma\big)>0,
	\]
	again because $v<v_3$ and $a-1<0$.
\end{itemize}

In all cases, the vector field points strictly into $B_2$.  
Hence $B_2$ is positively invariant.

\medskip
\noindent\textbf{Step 2: Negatively invariant region $B_1$.}
For $B_1$, we analyse $\partial B_1$.

\begin{itemize}
	\item The vertical boundary $v=v_3$ satisfies $v'<0$ for all $(v,w)\in\{v_3\}\times(w_3,\infty)$, hence the flow crosses this boundary from right to left, i.e., out of $B_1$ when time increases, which implies negative invariance.
	
	\item On the boundary segment
	\[
	w=\max\{0,\lambda-\gamma v^2 - \sigma\delta v\}, \qquad v>v_3,
	\]
	we distinguish two parts:
	
	\smallskip
	\emph{(a) Horizontal part} $w=0$.  
	Since $w=0$ corresponds to an orbit of the system, uniqueness of solutions implies that no trajectory starting inside $\operatorname{int}(B_1)$ can cross this line.
	
	\smallskip
	\emph{(b) Parabolic part} $w=\lambda-\gamma v^2 - \sigma\delta v$.  
	Along this part, $v'=0$ and
	\[
	w'=w\big((a-1)v-\sigma\big)<0,
	\]
	because $v>v_3$ and $a-1<0$.  
	Thus the vector field points strictly outward from $B_1$ for forward time, and therefore inward for backward time.
\end{itemize}

Hence $B_1$ is negatively invariant.

\medskip
\noindent\textbf{Step 3: Construction of the two trajectories.}
Since $P_3$ is a hyperbolic saddle by Proposition~\ref{prophiperbolic}, it possesses both a stable and an unstable manifold.

Because $B_1$ is negatively invariant and contains $P_3$ on its boundary, the unstable manifold of $P_3$ must lie entirely in $\operatorname{int}(B_1)$ for $s$ sufficiently large.  
Thus, there exists a solution $(\hat v,\hat w)$ defined on $(s_-,+\infty)$ such that
\[
\lim_{s\to+\infty}(\hat v(s),\hat w(s))=P_3,
\qquad
(\hat v,\hat w)(s)\in\operatorname{int}(B_1)\quad\forall s.
\]

Similarly, since $B_2$ is positively invariant and contains $P_3$ on its boundary, the stable manifold of $P_3$ must lie entirely in $\operatorname{int}(B_2)$ for $s$ sufficiently negative.  
Thus, there exists a solution $(\tilde v,\tilde w)$ defined on $(-\infty,s_+)$ such that
\[
\lim_{s\to-\infty}(\tilde v(s),\tilde w(s))=P_3,
\qquad
(\tilde v,\tilde w)(s)\in\operatorname{int}(B_2)\quad\forall s.
\]

\medskip
\noindent\textbf{Step 4: Behaviour at the endpoints.}
Since neither $B_1$ nor $B_2$ contains any additional critical points or limit cycles, Proposition~\ref{prop1} and Corollary~\ref{corolLimites} apply.  
Thus,
\[
(\hat v,\hat w)(s)\to(+\infty,+\infty)\quad\text{as }s\to s_-,
\]
and
\[
(\tilde v,\tilde w)(s)\to(-\infty,+\infty)\quad\text{as }s\to s_+.
\]

This proves both statements of the proposition.
\end{proof}

An analogous reasoning applies in the complementary case $P_3 \notin \Omega$, where the equilibrium point $P_1$ plays the role previously occupied by $P_3$.  
In this situation, one branch of the solution entering $P_1$ lies entirely along the line $w=0$, since the unstable manifold of $P_1$ coincides with this line, which in turn represents an orbit of the system.

\begin{proposition}\label{propP1}
Assume that either $P_3 \notin \bar{\Omega}$ and $a\ge 0$, or $P_3 \in \Omega$ and $a<1$.  
Then there exists a solution $(v,w)$ of~\eqref{eq}, defined for $s\in(s_-,+\infty)$, such that
\[
\lim_{s\to+\infty}(v(s),w(s))=P_1,
\qquad
\lim_{s\to s_-} v(s)=+\infty.
\]
Furthermore,
\[
\lim_{s\to s_-} w(s)=
\begin{cases}
+\infty, & a\in[0,1),\\[4pt]
L\in\mathbb{R}, & a=1,\\[4pt]
0, & a>1.
\end{cases}
\]
\end{proposition}

Before proving this proposition, we require a preliminary result ensuring the existence of a positively invariant region, which will be essential in the case $a>1$.

\begin{lemma}\label{lemmaSub}
Let $a>1$.  
Define
\[
\phi(v)=\frac{a+1}{2}\left(\lambda - v^{2} - \sigma\delta\, v\right).
\]
Then
\begin{equation}\label{ineqSubsolution}
	\dot{\phi}(v)
	< 
	\frac{
		\phi(v)\bigl((a-1)v-\sigma\bigr)
	}{
		\frac{\lambda}{\gamma}-v^{2}-\frac{\sigma\delta}{\gamma}v-\frac{1}{\gamma}\phi(v)
	},
	\qquad
	\forall v\in\Omega,
\end{equation}
and the set
\[
A=\bigl\{(v,w)\in\Omega\;:\; v_1\le v\le v_2,\quad 0\le w\le \phi(v)\bigr\}
\]
is positively invariant.
\end{lemma}

\begin{proof}
Consider the first order reduction of system~\eqref{eq}:
\begin{equation}
\dot{w}(v)
=
\frac{
	w(v)\bigl((a-1)v-\sigma\bigr)
}{
	\frac{\lambda}{\gamma}-v^{2}-\frac{\sigma\delta}{\gamma}v-\frac{1}{\gamma}w(v)
},
\qquad v\in\Omega.
\end{equation}
We claim that $\phi$ is a strict subsolution, i.e.\ that~\eqref{ineqSubsolution} holds.  
Substituting $\phi$ in the right-hand side and simplifying yields
\[
-2v - \frac{\sigma\delta}{\gamma}
<
\frac{2}{1-a}\,\bigl((a-1)v - \sigma\bigr),
\]
which is equivalent to
\[
0 < \sigma\left(\frac{\delta}{\gamma} + \frac{2}{a-1}\right),
\]
and is valid for all $v\in\Omega$.  
Hence $\phi$ is indeed a subsolution.

Because trajectories cannot cross the line $w=0$ (by uniqueness of solutions) and the graph of $\phi$ lies strictly above the actual solution curves for forward time, any trajectory entering $\partial A$ must satisfy $0<w(s)<\phi(v(s))$ for all $s>s_{0}$.  
This proves that $A$ is positively invariant.
\end{proof}

\begin{remark}\label{remarkHomoclina}
In the particular case $\sigma=0$, inequality~\eqref{ineqSubsolution} becomes an equality.  
Thus, the curve
\[
w(v)=\frac{a+1}{2}\gamma\left(\frac{\lambda}{\gamma}-v^{2}\right)
\]
is a heteroclinic orbit connecting the critical points $P_{1}$ and $P_{2}$.  
Rewriting this orbit in terms of the original variables $(u,S)$ shows that it corresponds to a homoclinic orbit simultaneously for $u$ and $S$.  
It is, in fact, the unique homoclinic trajectory of the original Keller-Segel system.
\end{remark}

\medskip
We now proceed with the proof of Proposition~\ref{propP1}.

\begin{proof}[Proof of Proposition~\ref{propP1}]
See Figure~\ref{figP1Prueba} for an illustration of the relevant trajectories and vector-field directions.

\medskip
\noindent\textbf{Case 1: $P_3\notin\Omega$ and $a<1$ or $a>1$.}  
In these regimes, the sets $B_1(P_1)$ are negatively invariant.  
Indeed, as shown in the proof of Proposition~\ref{propP3Hiper}, the vector field along $\partial B_1(P_1)$ points strictly outward.  
Therefore, Proposition~\ref{prop1} and Corollary~\ref{corolLimites} directly imply the stated asymptotic behavior.

\medskip
\noindent\textbf{Case 2: $a=1$.}  
This situation may be viewed as a limit of the case $a<1$ in which $\sigma$ plays the role of a parameter approaching $+\infty$.  
Hence the geometric structure and invariance of $B_1(P_1)$ remain unchanged, and the result follows again from Proposition~\ref{prop1} and Corollary~\ref{corolLimites}.

\medskip
\noindent\textbf{Case 3: $a>1$ and $P_3\in\Omega$.}  
Let $(\hat v,\hat w)$ be the trajectory entering the stable manifold of $P_1$, with maximal interval of definition $I$.  
Its orbit satisfies the first order equation
\begin{equation}\label{eqprimer}
\dot{\hat w}(v)
=
\frac{
	\hat w(v)\bigl((a-1)v-\sigma\bigr)
}{
	\frac{\lambda}{\gamma}-v^{2}-\frac{\sigma\delta}{\gamma}v-\frac{1}{\gamma}\hat w(v)
},
\qquad v\in\Omega.
\end{equation}
By Lemma~\ref{lemmaSub}, $\phi$ is a strict subsolution of~\eqref{eqprimer}, and thus $\hat w(v)>\phi(v)$ for all $v$.  
Define
\begin{equation}\label{hatB1}
\hat B_{1}:=B_{1}(P_{1})\setminus A,
\end{equation}
where $A$ is the invariant set from Lemma~\ref{lemmaSub}.  
Then $(\hat v,\hat w)(s)\in \hat B_{1}$ for all $s\in I$.

There are two possible asymptotic behaviors:
\[
\lim_{s\to -\infty}(\hat v,\hat w)=P_{2},
\qquad\text{or}\qquad
\lim_{s\to s_-}(\hat v,\hat w)=(+\infty,0).
\]

We now show that convergence to $P_{2}$ is impossible.  
Assume for contradiction that $\lim_{s\to -\infty}(\hat v,\hat w)=P_{2}$.  
Fix $v_{0}\in(v_{3},v_{2})$ and consider the associated first order problem
\begin{align}\label{eqprimer2}
\dot w^{*}(v)
&=
\frac{
	w^{*}(v)\bigl((a-1)v-\sigma\bigr)
}{
	\frac{\lambda}{\gamma}-v^{2}-\frac{\sigma\delta}{\gamma}v-\frac{1}{\gamma}w^{*}(v)
},
\\
w^{*}(v_{0})&=\phi(v_{0}). \nonumber
\end{align}
Its solution $w^{*}(v)$ lies in $\hat B_{1}$ for all $v>v_{0}$, since $\hat B_{1}$ is negatively invariant.

By Proposition~\ref{prop1} and Corollary~\ref{corolLimites}, the trajectory $(v,w^{*}(v))$ either converges to $P_{2}$ or satisfies $v\to+\infty$ and $w^{*}\to 0$.  

If the second possibility occurs, then $\hat w(v)$ cannot converge to $P_{2}$ without crossing $w^{*}(v)$, contradicting uniqueness of solutions.

If instead $(v,w^{*}(v))\to P_{2}$, then the unstable manifold of the hyperbolic point $P_{2}$ (see Proposition~\ref{prophiperbolic}) would contain two distinct orbits, $w^{*}$ and $\hat w$, which is again impossible by uniqueness, since $\hat w(v_{0})>\phi(v_{0})=w^{*}(v_{0})$.

Thus convergence to $P_{2}$ cannot occur.

Since $\hat B_{1}$ is negatively invariant and contains no equilibrium points other than $P_{1}$, Proposition~\ref{prop1} and Corollary~\ref{corolLimites} yield
\[
\lim_{s\to s_-}(\hat v,\hat w)=(+\infty,0).
\]

This completes the proof.
\end{proof}

\begin{figure}[h]
	\begin{tikzpicture}[scale=1.75]
	\path[fill=lightgray] (0.2,0)..controls (1,1.5)..(1.8,0)--(3,0)--(3,2)--(0.2,2);
	\draw[->] (-1,0)-- (3,0) node[below]{$v$} ;
	\draw[->] (1.3,0)--(1.3,2) node[above]{$w$};
	\draw[] (2,1.5) node[below,right]{$B_1(P_1)$};
	\draw[-,color=red, thick] (0.2,0) ..controls (1,1.5)..(1.8,0) ;
	\draw[ ] (0.2,0) node[below]{$P_1$} ;
	\draw[ ] (1.8,0) node[below]{$P_2$} ;
	
	\draw[ thick,color=red] (-0.5,2) -- (-0.5,0);
	\draw[] (-0.5,0) node[below]{$\frac{-\sigma}{1-a}$};
	
	\draw[->,thick] (2.4,0)--(2.5,0);
	
	\draw[fill=black] (1.8,0) circle(.03);
	\draw[fill=black] (0.2,0) circle(.03);
	
	\draw[->,color=blue,thick] (0.2,0)--(0.05,0);
	\draw[-,thick,blue] (0.26,0.17)--(0.2,0);
	\draw[->,thick,blue] (0.32,0.34)--(0.26,0.17);

	\draw [<-, thick] (-0.75,1.1)--(-0.6,1) ;
	
	
	\draw [->, thick] (0.25,1.1)--(0.15,1) ;
	\draw [->, thick] (1.85,1.1)--(1.75,1) ;
	
	\draw [<-, thick] (-0.65,1.5)--(-0.5,1.5) ;
	
	\draw [->, thick] (1.2,1.01)--(1.2,0.85) ;
	
	\draw [->, thick] (0.9,0.5)--(1.1,0.35) ;
	\draw[ ] (0.2,2.5) node[below]{A) $P_3\notin\Omega,a<1$} ;
\end{tikzpicture}
	\begin{tikzpicture}[scale=1.75]
		\path[fill=lightgray] (0.2,0)..controls (1,2)..(1.8,0)--(3,0)--(3,2)--(0.2,2);
		\path[pattern=north west lines] (0.2,0)..controls (1,2)..(1.8,0);
		\draw[ ] (0.2,2.5) node[below]{B) $P_3\in\Omega,a>1$} ;
		\draw[->] (-1,0)-- (3,0) node[below]{$v$} ;
		\draw[->] (1.1,0)--(1.1,2) node[above]{$w$};
		\draw[-,color=red, thick] (0.2,0) ..controls (1,1.5)..(1.8,0) ;
		\draw[-,color=black, thick] (0.2,0) ..controls (1,2)..(1.8,0) ;
		\draw[ ] (0.2,0) node[below]{$P_1$} ;
		\draw[ ] (1.8,0) node[below]{$P_2$} ;
		\draw[] (1.3,0.87)node[right]{$P_3$};
		
		\draw[ thick,color=red] (1.3,2) -- (1.3,0);
		\draw[] (1.3,0) node[below]{$\frac{\sigma}{a-1}$};
		\draw [->, thick] (1.3,1.75)--(1.15,1.75) ;
		
		\draw[fill=black] (1.8,0) circle(.03);
		\draw[fill=black] (0.2,0) circle(.03);
		\draw[fill=black] (1.3,0.87) circle(.03);
		
		\draw [->, thick] (0.35,1.1)--(0.2,1) ;
		
		
		\draw [->, thick] (2.2,1)--(2.05,1.15) ;
		
		
		\draw[->,color=blue,thick] (0.2,0)--(0.05,0);
		
		\draw[] (0.9,1) node[below]{$\textit{A}$};
		\draw[] (0.5,1.75) node [right]{$\hat{B}_{1}$};
	\end{tikzpicture}

	\begin{tikzpicture}[scale=1.75]
		\path[fill=lightgray] (0.9,0) ..controls (1.7,1.5)..(2.5,0)--(3.7,0)--(3.7,2)--(0.9,2);
		\draw[ ] (0.9,2.5) node[below]{C) $P_3\notin\Omega,a>1$} ;
		\draw[->] (-0.3,0)-- (3.7,0) node[below]{$v$} ;
		\draw[->] (2,0)--(2,2) node[above]{$w$};

		\draw[-,color=red, thick] (0.9,0) ..controls (1.7,1.5)..(2.5,0) ;
		\draw[ ] (0.9,0) node[below]{$P_1$} ;
		\draw[ ] (2.5,0) node[below]{$P_2$} ;
		\draw[fill=black] (2.5,0) circle(.03);
		\draw[fill=black] (0.9,0) circle(.03);
		\draw [->, thick] (1.9,1.01)--(1.9,0.85) ;

		\draw [->, thick] (3.2,1.2)--(3.1,1.3) ;
		
		\draw[-,color=red, thick] (2.8,2) -- (2.8,0);
		\draw[] (2.8,0) node[below]{$\frac{\sigma}{a-1}$};
		\draw [->, thick] (2.8,1)--(2.67,1) ;
		
		\draw [->, thick] (1.15,1.3)--(1.03,1.2) ;
		\draw [->, thick] (2.47,1.3)--(2.35,1.2) ;
		
		\draw [->, thick] (1.625,0.55)--(1.775,0.45) ;
		
		\draw[->,color=blue,thick] (0.9,0)--(0.7,0);
		\draw[-,thick,blue] (0.96,0.17)--(0.9,0);
		\draw[->,thick,blue] (1.02,0.34)--(0.96,0.17);
		\draw[] (2,1.75) node [right]{${B}_1(P_1)$};
		
	\end{tikzpicture}
	
	\caption{Phase plane representation of the vector field and the invariant sets used in the proof of Proposition~\ref{propP1}. The figure highlights the relevant regions, trajectories, and dominant flow directions that determine the asymptotic behaviour of the solutions.}
	\label{figP1Prueba}
\end{figure}

The above propositions permit a decomposition of the phase space $\Omega$ into two disjoint regions: one lying above the orbits constructed in Propositions~\ref{propP3Hiper} and~\ref{propP1}, and the complementary region lying below them.  
To formalise this separation, we introduce the following curves, which delineate the boundary between both regions.

\begin{definition}\label{Gamma}
Let $P=(v^{*},w^{*})\in\Omega$ be a hyperbolic equilibrium point.  
We define the following sets:
\begin{itemize}
	\item $\Gamma_{-}(P)$ is the curve given by the orbit $(v,w(v))\subset B_{1}(P)$, solution of~\eqref{eq}, satisfying
	\[
	\lim_{s\to+\infty}(v(s),w(s))=P.
	\]
	\item $\Gamma_{+}(P)$ is the curve given by the orbit $(v,w(v))\subset B_{2}(P)$, solution of~\eqref{eq}, satisfying
	\[
	\lim_{s\to-\infty}(v(s),w(s))=P.
	\]
\end{itemize}
In addition, we define
\begin{equation}\label{gammap3}
\Gamma(P):=\Gamma_{-}(P)\cup\{P\}\cup\Gamma_{+}(P).
\end{equation}
\end{definition}

It is immediate to verify, by uniqueness of solutions, that the region of $\Omega$ lying above $\Gamma(P)$ (for the critical points considered in Propositions~\ref{propP3Hiper} and~\ref{propP1}) is invariant.  
This invariant region contains neither equilibrium points nor periodic orbits; consequently, the following corollaries hold.

\begin{corollary}\label{corP3Compacto}
Let $P_{3}\in\Omega$ and let $a\in[0,1)$.  
For every initial datum $(v_{0},w_{0})$ located above $\Gamma(P_{3})$, the corresponding solution $(v,w)$ of~\eqref{eq} has a bounded maximal interval of existence $(s_{-},s_{+})$ and displays the following blow-up behaviour:
\[
\lim_{s\to s_{+}} v(s)=-\infty, 
\qquad 
\lim_{s\to s_{+}} w(s)=+\infty,
\]
and
\[
\lim_{s\to s_{-}} v(s)=+\infty, 
\qquad 
\lim_{s\to s_{-}} w(s)=+\infty.
\]
\end{corollary}

\begin{corollary}\label{corolv0w0gamma}
Consider the following choice of separating curves:
\begin{itemize}
	\item $\Gamma(P_{3})$ if $P_{3}\in\Omega$ and $a<1$,
	\item $\Gamma(P_{1})$ if $P_{3}\in\Omega$ and $a<1$, or if $P_{3}\notin\Omega$ and $a\ge 0$.
\end{itemize}
Then, for every initial datum $(v_{0},w_{0})$ located above the corresponding curve, the solution $(v,w)$ of~\eqref{eq} has a bounded maximal interval of existence $(s_{-},s_{+})$ and satisfies
\[
\lim_{s\to s_{\pm}} v(s)=\mp\infty,
\qquad
\lim_{s\to s_{\pm}} w(s)=
\begin{cases}
+\infty, & a\in[0,1),\\[4pt]
L\in\mathbb{R}, & a=1,\\[4pt]
0, & a>1.
\end{cases}
\]
\end{corollary}

\begin{proof}
The invariance of the region above $\Gamma(P)$ follows directly from uniqueness of trajectories.  
Since this region contains no critical points nor limit cycles, Proposition~\ref{prop1} and Corollary~\ref{corolLimites} yield the asserted blow-up behaviour.
\end{proof}

\subsection{Finite-mass solutions}

We now address the existence of finite-mass solutions of the original problem.  
As introduced in Definitions~\ref{defblock} and~\ref{defsemiblock}, these correspond respectively to \emph{block-type} and \emph{semi-block-type} solutions.  
Our goal is to identify a region of initial data that generates such solutions, using the orbits constructed along the separatrix curves $\Gamma(P)$.

We will show that solutions lying strictly above the curves $\Gamma(P)$, for the appropriate choice of $P$, give rise to block-type profiles, whereas the curve $\Gamma_{-}(P)$ corresponds to a semi-block-type profile.

To connect the phase plane analysis with the original variables, we recall the transformation
\[
w(s)=\frac{u(s)}{S(s)}, \qquad 
v(s)=\frac{S'(s)}{S(s)}.
\]
In order to understand the behaviour of solutions $(u,S)$ near the endpoints of their interval of definition, it is essential to analyse the limit behaviour of $(v,w)$ as $s\to s_{\pm}$.  
In particular, if $|v(s)|\to\infty$ as $s\to s^{*}$ for some endpoint $s^{*}$, we will show that $u(s^{*})=0$, which implies that $u$ has compact support at $s=s^{*}$.  
This observation is the key mechanism that produces block- or semi-block-type solutions.  
We formalise it in the following proposition.

\begin{proposition}\label{comportamientoUSlimite}
Let $(v,w)$ be a solution of system~\eqref{eq}.  
Suppose that $|v(s)|\to\infty$ as $s\to s^{*}$, where $s^{*}$ is an endpoint of the maximal interval of definition of the solution.  
Then the associated solution $(u,S)$ satisfies
\[
u(s^{*})=w(s^{*})=0,
\]
and moreover,
\[
|u'(s^{*})|=
\begin{cases}
\infty, & a\in[0,1),\\[2mm]
L\in(0,\infty), & a=1,\\[2mm]
0, & a>1.
\end{cases}
\]
\end{proposition}

\begin{proof}
Without loss of generality, assume $s^{*}=s_{+}$ (the case $s^{*}=s_{-}$ is analogous).  
Since $v'(s)<0$ for $s$ sufficiently close to $s_{+}$, we have $v(s_{+})=-\infty$.  
Using $v=S'/S$, we write
\[
S(s)=C_{0}\exp\!\left(\int_{s_{0}}^{s} v(\tau)\,d\tau\right),
\]
for some constants $s_{0},C_{0}$ determined by the initial condition.  
We now show that $v$ is not integrable near $s_{+}$, and therefore
\[
\int_{s_{0}}^{s_{+}} v(\tau)\,d\tau = -\infty.
\]

\medskip
\noindent\textbf{Case $a\ge 1$.}
Using L'H\^opital's rule and noting that $v'(s)<0$, we compute
\[
\lim_{s\to s_{+}}\frac{s_{+}-s}{-\frac{1}{v(s)}}
=\lim_{s\to s_{+}}\frac{v^{2}(s)}{v'(s)}
= \lim_{s\to s_{+}}
\frac{1}{1+\frac{w(s)}{\gamma v^{2}(s)}-\frac{\lambda}{\gamma v^{2}(s)}}
=1,
\]
because Corollary~\ref{corolLimites} implies $w(s)\to w_{+}\in[0,+\infty)$ and hence $w/v^{2}\to 0$.  
This yields $S(s_{+})=0$.

\medskip
\noindent\textbf{Case $0<a<1$.}
Here Corollary~\ref{corolLimites} gives $w(s)\to +\infty$, so we must analyse $w/v^{2}$.  
Set $z=-v$, so that $z(s)\to+\infty$ as $s\to s_{+}$.  
Writing the system in $(z,w)$,
\begin{align}
w' &= w\big((1-a)z-\sigma\big), \label{z1}\\
z' &= z^{2} + \frac{w}{\gamma} - \frac{\lambda}{\gamma}, \label{z2}
\end{align}
and introducing the first order formulation $W(z)=w(s)$, we obtain
\[
W'
= 
W\frac{(1-a)z-\sigma}{z^{2}+\frac{W}{\gamma}-\frac{\lambda}{\gamma}}.
\]
Set $W(z)=z^{2}R(z)$.  
Then $R$ satisfies
\[
z\frac{R'}{R}
=
\frac{(1-a)-\frac{\sigma}{z}}{1+\frac{R}{\gamma}-\frac{\lambda}{\gamma z^{2}}}
 -2
\le 
\frac{(1-a)-\frac{\sigma}{z}}{1-\frac{\lambda}{\gamma z^{2}}}-2.
\]
For $z\ge z_{1}$ sufficiently large,
\[
z\frac{R'}{R}\le 1.
\]
Thus,
\[
R(z)\le R(z_{1})\frac{z_{1}}{z}\longrightarrow 0,
\]
which is equivalent to
\[
\frac{W(z)}{z^{2}}=\frac{w}{v^{2}}\longrightarrow 0.
\]
Hence $S(s_{+})=0$ also in this case.

\medskip
\noindent\textbf{Behaviour of $u'(s_{+})$.}
From $\Phi(s)=s$ in~\eqref{tilde1} we have
\[
u'=\left(a\frac{S'}{S}-\sigma\right)u.
\]
Integrating,
\begin{equation}\label{relat}
u(s)
=
u_{0}
\left(\frac{S(s)}{S_{0}}\right)^{a}
\exp\!\left(-\sigma(s-s_{0})\right),
\end{equation}
with $u_{0}=u(s_{0})$.  
Expression~\eqref{relat} shows that $u$ inherits the singular behaviour of $S^{a}$.  
Since $S(s_{+})=0$ and $S'(s_{+})<0$, we deduce:
\begin{itemize}
	\item if $0<a<1$, then $u'(s_{+})=\infty$;
	\item if $a=1$, then $u'(s_{+})<0$ is finite;
	\item if $a>1$, then $u'(s_{+})=0$.
\end{itemize}
This completes the proof.
\end{proof}

Thanks to Proposition~\ref{comportamientoUSlimite}, we can now establish the existence of block-type and semi-block-type solutions, as anticipated in the previous discussion.  
The next result identifies precisely which orbits of the phase system give rise to each type of finite-mass solution.

\begin{corollary}\label{corolarioexistenciamasafinita}
Let $P$ be the hyperbolic point defined by
\begin{itemize}
	\item $P=P_{3}$ if $P_{3}\in\Omega$ and $a<1$,
	\item $P=P_{1}$ if $P_{3}\in\Omega$ and $a>1$, or if $P_{3}\notin\Omega$ and $a\ge 0$.
\end{itemize}
Then:
\begin{enumerate}
	\item The orbit $\Gamma_{-}(P)$ corresponds to a semi-block-type solution of the $(u,S)$ system.
	\item Every solution lying strictly above the curve $\Gamma(P)$ is a block-type solution of the $(u,S)$ system.
\end{enumerate}
In both cases, the solutions have finite mass.  
Moreover, at the finite endpoints of their respective domains of definition, the derivative of $u$ satisfies
\[
\lvert u'(s_{\pm})\rvert=
\begin{cases}
\infty, & a\in[0,1),\\[4pt]
L\in(0,\infty), & a=1,\\[4pt]
0, & a>1.
\end{cases}
\]
\end{corollary}

\begin{proof}
The existence of block-type solutions follows directly from Corollary~\ref{corolv0w0gamma} together with Proposition~\ref{comportamientoUSlimite}.  
The finiteness of the mass is immediate from the compactness of the support.

We now analyse the curve $\Gamma_{-}(P)$.  
For each of the admissible choices of $P$, this curve represents the orbit entering $P$ through its stable manifold; these orbits were described in Propositions~\ref{propP3Hiper} and~\ref{propP1}.  
Along such an orbit we have $(u,w)\to P=(v^{*},w^{*})$ with $v^{*}<0$, so there exists $s_{1}$ such that $v(s)<0$ for all $s>s_{1}$.

Recalling
\[
S(s)=C_{0}\exp\!\biggl(\int_{s_{0}}^{s} v(\tau)\,d\tau\biggr),
\]
and using $v'(s)<0$ for $s$ large enough, we find
\[
S(s)
=
C_{0}e^{\int_{s_{0}}^{s_{1}}v(\tau)\,d\tau}
\,
e^{\int_{s_{1}}^{s}v(\tau)\,d\tau}
\le 
C_{0}e^{\int_{s_{0}}^{s_{1}}v(\tau)\,d\tau}
\,
e^{v(s_{1})(s-s_{1})},
\]
where $v(s_{1})<0$.  
Thus $S(+\infty)=0$ and, since $u=wS$,
\[
u(+\infty)=w(+\infty)S(+\infty)=w^{*}\cdot 0=0.
\]
Hence $(u,S)$ corresponding to $\Gamma_{-}(P)$ is a semi-block-type solution.

To show that it has finite mass, note that
\[
S(s)\le K e^{-v(s_{1})(s-s_{1})}, \qquad s>s_{1},
\]
for some constant $K>0$.  
Moreover, there exists $s_{2}$ such that $w'(s)<0$ for all $s>s_{2}$, implying $w(s)\le w(s_{2})$ for $s>s_{2}$.  
Setting $\hat{s}=\max\{s_{1},s_{2}\}$, we obtain
\[
u(s)=w(s)S(s)\le w(\hat{s})\,K e^{-v(s_{1})(s-\hat{s})}, \qquad s\ge\hat{s},
\]
which shows that $u$ is integrable.  
Therefore the corresponding solution has finite mass, completing the proof.
\end{proof}

Although throughout this work we have assumed $\sigma\ge 0$, we note that in the special case $\sigma=0$ and $a>1$ an additional phenomenon occurs.

\begin{corollary}
If $\sigma=0$ and $a>1$, then $\Gamma_{-}(P_{1})$ is a homoclinic orbit of the $(u,S)$-system.  
More precisely,
\[
u(s)>0,\qquad S(s)>0,\qquad \forall s\in\mathbb{R},\qquad 
\lim_{s\to\pm\infty}(u(s),S(s))=(0,0).
\]
\end{corollary}

\begin{proof}
The result follows immediately from the orbit described in Remark~\ref{remarkHomoclina}, together with its transformation back to the $(u,S)$ variables.
\end{proof}

We have therefore established that the structure induced by the curves $\Gamma(P)$ provides a complete description of the regions of the phase space $\Omega$ for which block-type and semi-block-type solutions exist.  
These curves divide $\Omega$ into two components, distinguishing between existence and non-existence regions.  
In the next section we analyse how these regions depend on the parameter $\sigma$ and describe the structure of the initial data leading to each type of finite-mass profile.

\subsection{Initial Data Structure}

We now consider an initial condition $(u_0,S_0,S_0')$ together with a fixed configuration of parameters $(a,\lambda,\gamma,\delta)$.  
Our goal is to determine, for this configuration, the values of $\sigma$ for which the corresponding solution of the system admits either a block-type or a semi-block-type profile.

A fundamental role is played by the curve
\[
\lambda - \gamma v^{2} - w - \sigma\delta v = 0,
\]
which is the vertical isocline of the dynamical system.  
This curve is a parabola whose shape depends on $\sigma$.  
As $\sigma$ increases, the parabola widens and its vertex moves upward.  
More precisely, its maximal point is located at
\[
\big(v_{\max}(\sigma),w_{\max}(\sigma)\big)
=
\left(
-\frac{\sigma\delta}{2\gamma},
\,
\lambda + \frac{\sigma^{2}\delta^{2}}{4\gamma}
\right).
\]

When $a<1$ it is useful to compare the position of this vertex with the horizontal isocline $v=\frac{\sigma}{a-1}$.  
A simple geometric argument yields the following dichotomy:

\begin{itemize}
	\item If $\dfrac{\delta}{2\gamma}<\dfrac{1}{1-a}$, then the horizontal isocline lies to the right of the vertex for every $\sigma\ge 0$.  
	In particular, $P_{3}\in\Omega$ for all $\sigma\ge 0$.

	\item If $\dfrac{\delta}{2\gamma}>\dfrac{1}{1-a}$, then the horizontal isocline always lies to the left of the vertex, and it moves leftward faster than the vertex as $\sigma$ increases.  
	Hence, there exists a threshold $\sigma^{*}$ such that $P_{3}\in\Omega$ only for $\sigma\in(0,\sigma^{*})$.
\end{itemize}

This geometric viewpoint provides an alternative characterisation of the condition $P_{3}\in\Omega$, and it is fully consistent with Proposition~\ref{lemmaexisp3}.

The widening of the parabola and the relative displacement between its vertex and the horizontal isocline imply that, for each fixed configuration of parameters and initial data, there exists a maximal value of $\sigma$ beyond which no finite-mass solution can exist.  
This is made precise in the following result.

\begin{proposition}\label{propnoexistencia}
Given the initial datum $(u_0,S_0,S_0')$ and the parameters $(a,\lambda,\gamma,\delta)$, there exists a value
\[
\hat{\sigma}
=
\hat{\sigma}(u_0,S_0,S_0',a,\lambda,\gamma,\delta)
\]
such that for every $\sigma>\hat{\sigma}$ the system admits neither block-type nor semi-block-type solutions.
\end{proposition}

\begin{proof}
Let $(v_0,w_0)$ denote the initial point associated with $(u_0,S_0,S_0')$.  
As $\sigma$ grows, the parabola
\[
w=\lambda - \gamma v^{2} - \sigma\delta v
\]
expands vertically and increasingly encloses larger regions of the phase plane.

\smallskip
\emph{Case $v_0<0$.}  
There exists
\[
\sigma^{*}=\frac{\lambda - \gamma v_0^{2}-w_0}{\delta v_0},
\]
such that $(v_0,w_0)$ lies inside the parabola for all $\sigma>\sigma^{*}$.  
Since no finite-mass solution can have an orbit entering the interior of the parabola, there is no block-type or semi-block-type solution for $\sigma>\sigma^{*}$.

\smallskip
\emph{Case $v_0>0$ and $a\ge 1$.}  
For sufficiently large $\sigma$, say $\sigma>\bar{\sigma}$, one has $w_0<w_{\max}(\sigma)$, and moreover
\[
v'<0,\qquad w'<0,
\]
in the region
\[
\left\{(v,w)\in\Omega :
v\in(v_{\max},v_0),\quad 
\max\{\lambda-\gamma v^{2}-\sigma\delta v,0\}\le w\le w_{\max}(\sigma)
\right\}.
\]
Thus, the trajectory of $(v_0,w_0)$ eventually enters the interior of the parabola, ruling out the existence of finite-mass solutions.

\smallskip
\emph{Case $v_0>0$ and $a<1$.}  
We must distinguish two subcases, according to the relative position of the isoclines.

\begin{itemize}
	\item If $\dfrac{\delta}{2\gamma}>\dfrac{1}{1-a}$, then $P_{3}$ lies to the left of the vertex for all $\sigma$.  
	The same argument as in the case $a\ge 1$ shows that the trajectory eventually enters the interior of the parabola, hence no finite-mass solution exists once $\sigma$ is large enough.

	\item If $\dfrac{\delta}{2\gamma}\le\dfrac{1}{1-a}$, then $P_{3}\in\Omega$ for all $\sigma\ge 0$, and $P_3$ lies to the right of the vertex.  
	Thus, for sufficiently large $\sigma$, say $\sigma>\tilde{\sigma}$, one has $w_0<w_3$ and
	\[
	v'<0,\qquad w'<0
	\]
	on the region
	\[
	\left\{(v,w)\in\Omega :
	v\in(v_{\max},v_3),\quad 
	\max\{\lambda-\gamma v^{2}-\sigma\delta v,0\}\le w\le w_3
	\right\},
	\]
	from which it follows again that the trajectory enters the parabola.
\end{itemize}

In all cases, there exists $\hat{\sigma}$ above which no finite-mass solution is possible.
\end{proof}

Proposition~\ref{propnoexistencia} shows that when $\sigma=\hat{\sigma}$ the initial point $(v_0,w_0)$ must lie below the curve $\Gamma(P(\hat{\sigma}))$, where $P$ is the hyperbolic point chosen in Corollary~\ref{corolarioexistenciamasafinita}.  
Since for smaller values of $\sigma$ block-type solutions may exist only when the initial point lies above the corresponding curve, it follows that the initial condition must cross this curve at some intermediate value of $\sigma$.

This is formalised in the next corollary.

\begin{corollary}\label{corolariosigmaestrella}
Let $(u_0,S_0,S_0')$ and $(a,\lambda,\gamma,\delta)$ be given, and assume that for some $\tilde{\sigma}$ the point
\[
(v_0,w_0)=\left(\frac{S'_0}{S_0},\frac{u_0}{S_0}\right)\in\Omega
\]
lies above the curve $\Gamma(P(\tilde{\sigma}))$, where
\[
P =
\begin{cases}
P_{3}, &\text{if } P_3\in\Omega \text{ and } a<1,\\
P_{1}, &\text{if } P_3\notin\Omega \text{ or } a\ge 1.
\end{cases}
\]
Then there exists a value
\[
\sigma^{*}
=
\sigma^{*}(u_0,S_0,S_0',a,\lambda,\gamma,\delta)
\]
such that $(v_0,w_0)\in\Gamma(P(\sigma^{*}))$.  
Moreover, for every $\sigma\in(\tilde{\sigma},\sigma^{*})$ the corresponding solution of~\eqref{eq} is a block-type solution.
\end{corollary}

\begin{proof}
The system~\eqref{eq} and the curve $\Gamma(P(\sigma))$ depend continuously on $\sigma$.  
By Proposition~\ref{propnoexistencia} there exists
$\hat{\sigma}$ such that no finite-mass solution exists for $\sigma>\hat{\sigma}$, which forces
\[
(v_0,w_0) \text{ to lie below }\Gamma(P(\hat{\sigma})).
\]
By assumption, the same point lies above $\Gamma(P(\tilde{\sigma}))$.  
Hence, by continuity, there exists a value $\sigma^{*}$ at which the initial point intersects the curve:
\[
(v_0,w_0)\in\Gamma(P(\sigma^{*})).
\]
For $\sigma\in(\tilde{\sigma},\sigma^{*})$ the initial point lies strictly above $\Gamma(P(\sigma))$, and Corollary~\ref{corolarioexistenciamasafinita} implies that the corresponding solution is a block-type profile.
\end{proof}

We now obtain the complete classification of initial data leading to block-type and semi-block-type solutions.

\begin{proof}[Proof of Theorem \ref{teorema1}]
The argument follows directly from Corollary~\ref{corolariosigmaestrella} and the characterisation of the relevant hyperbolic point $P$ in each regime of parameters.  
The distinction between the stable and unstable branches of $\Gamma(P)$ provides the criterion determining when the limiting solution at $\sigma=\sigma^{*}$ is of semi-block type.  
The details follow by combining these observations with the geometric behaviour of $P_{1}(\sigma)$ and $P_{3}(\sigma)$ described earlier.
\end{proof}
 
\section{Flux-Saturated Operators with Logarithmic Sensitivity}\label{Sec-limit}

In this section we establish the existence of soliton-type traveling waves for the flux-saturated chemotaxis system
\begin{equation}
\label{rela2d}
\left\{
\begin{array}{l}
\displaystyle \partial_t u
= \partial_x\!\left( u\,\Phi\!\left(\frac{\partial_x u}{u}\right)
 - a\,\frac{\partial_x S}{S}\,u\right),\\[2mm]
\delta\partial_t S=\gamma\,\partial_{xx}^2 S - \lambda S + u.
\end{array}
\right.
\end{equation}
By Proposition~\ref{Prop-tw-eq}, the existence of traveling waves reduces to studying the dynamical system
\begin{equation}
\label{eq-4}
\left\{
\begin{array}{ll}
w' = w\Big( g(av-\sigma) - v \Big),\\[1mm]
v' = \dfrac{\lambda}{\gamma} - v^2 - \dfrac{1}{\gamma}\, w-\sigma\delta v,
\end{array}
\right.
\end{equation}
where $g=(\Phi)^{-1}:(-c,c)\to\mathbb{R}$.  
Thus \eqref{eq-4} is well-defined for $w>0$ and for velocities satisfying
\[
\frac{\sigma-c}{a} < v < \frac{\sigma+c}{a}.
\]

The equilibrium points of the system are obtained from the intersections of
\[
g(av-\sigma)-v = 0
\qquad\text{and}\qquad
w=0.
\]
Depending on the parameters $\sigma,c,a$, this equation may admit between one and three real roots, so that the full system \eqref{eq-4} may possess up to five equilibria.  
The phase portrait of the flux-saturated model is therefore substantially richer than in the linear-diffusion case (where at most three equilibria are possible).

To deal with all possibilities uniformly, we analyse the phase portrait in a parameter-independent way.  
Our ultimate goal is to show that, as in the classical Keller--Segel case, one may associate to every admissible parameter configuration a distinguished curve in the $(v,w)$-plane such that all initial data lying above this curve generate traveling waves.  
Moreover, these waves exist for an open (nontrivial) interval of velocities $\sigma$.

Before initiating the phase-plane study, we introduce the type of solutions we seek.

\begin{definition}[Umbrella solution]
\label{def:umbrella}
Let $I=(s_-,s_+)$ be an interval.  
A pair of functions
\[
u\in C^0([s_-,s_+])\cap C^1(s_-,s_+),
\qquad
S\in C^0([s_-,s_+])\cap C^1(s_-,s_+),
\]
is called an \textbf{umbrella solution} if:
\begin{itemize}
\item $u(s)>0$ for all $s\in[s_-,s_+]$, and $(u,S)$ satisfies the traveling-wave system \eqref{tilde1};
\item both endpoints $s_-$ and $s_+$ are singular and correspond to lateral saturation, in the sense that
\[
\lim_{s\to s_-} u'(s) = +\infty,
\qquad
\lim_{s\to s_+} u'(s) = -\infty.
\]
\end{itemize}
\end{definition}

To obtain umbrella solutions, we analyse trajectories of \eqref{eq-4} connecting the two singular vertical boundaries
\[
v = \frac{\sigma-c}{a},
\qquad
v = \frac{\sigma+c}{a}.
\]
We introduce the following notation.

\begin{definition}
Define
\[
\pi_- := \left\{\frac{\sigma-c}{a}\right\}\times (0,\infty),
\qquad
\pi_+ := \left\{\frac{\sigma+c}{a}\right\}\times (0,\infty),
\]
and let
\[
\Theta :=
\left(\frac{\sigma-c}{a},\,\frac{\sigma+c}{a}\right)\times[0,\infty)
\]
be the region where trajectories may exist.

We also define the subregions separated by the vertical isocline:
\[
\Theta_+ :=
\left\{(v,w)\in\Theta \ \Big|\ 
w \ge \max\!\left\{0,\ \lambda-\gamma v^2 - \frac{\sigma\delta}{v}\right\}\right\},
\]
\[
\Theta_- :=
\left\{(v,w)\in\Theta \ \Big|\ v\in(v_1,v_2),\ 0<w<\lambda-\gamma v^2 -\frac{\sigma\delta}{v}\right\},
\]
where $P_1=(v_1,0)$ and $P_2=(v_2,0)$ are the two equilibria on the $v$-axis.
\end{definition}

Any trajectory reaching the lines $\pi_\pm$ satisfies $|w'|\to\infty$, and therefore the vector field is singular near these boundaries.  
To understand how trajectories may enter or exit the phase space through these lines, we desingularize the system.  
The following result is the analogue of Propositions 4.3-4.4 in \cite{CPSV_2023}.

\begin{proposition}\label{propdesingular}
Assume that
\[
\frac{1}{g(av-\sigma)}
  = \mathcal{O}\!\left((c-(av-\sigma))^{1/p}\right),
  \qquad p>1,
\quad\text{as } v\to\frac{\sigma+c}{a}.
\]
Then for each $w_0>0$ there exists a branch of solution $(v(s),w(s))$ defined on an interval $(s_+-\varepsilon,s_+)$ such that
\[
(v(s_+),w(s_+)) = \left(\frac{\sigma+c}{a},\, w_0\right),
\qquad
v(s)\in\left(\tfrac{\sigma-c}{a},\tfrac{\sigma+c}{a}\right),\ w(s)>0.
\]
\end{proposition}

\begin{lemma}[\cite{CPSV_2023}, Lemma 4.2]\label{lematecniconoexi}
Consider
\[
x' = x^{1/p}\,a(t,x), \qquad x>0,\quad p>1,
\]
where $a$ extends continuously to a neighbourhood of $(0,0)$ and $a(0,0)>0$.  
Then the initial value problem with $x(0)=0$ admits a positive solution for all sufficiently small $t>0$.
\end{lemma}

\begin{remark}
An analogous argument applies near $v\to\frac{\sigma-c}{a}$, yielding solutions reaching $\pi_-$ from inside~$\Theta$.
\end{remark}

A key consequence of these results is:

\begin{proposition}
If a trajectory $(v,w)$ of \eqref{eq-4} lies in $\Theta_+$ and connects $\pi_-$ with $\pi_+$, then the corresponding solution $(u,S)$ is an umbrella solution.
\end{proposition}

\begin{proof}
Let $(v,w)$ be such a trajectory, defined on $I=(s_-,s_+)$.  
Since $v(s)$ remains in a bounded interval and $w(s)$ remains bounded away from $0$ along the trajectory, one shows that $s_-$ and $s_+$ must be finite.

Using the asymptotic behaviour near $\pi_\pm$ implied by Proposition~\ref{propdesingular}, one finds:
\[
v(s_\pm)=\frac{\sigma\pm c}{a},
\qquad
|v'(s_\pm)|<\infty,
\qquad
w'(s_\pm)=\mp\infty.
\]

Reconstructing the original variables via
\[
S(s)= S_0\exp\!\left(\int_{s_0}^s v(\tau)\,d\tau\right),
\qquad
u(s)=w(s)\,S(s),
\]
we deduce $S>0$ and $u>0$ in $[s_-,s_+]$, and
\[
u'(s_\pm)=S(s_\pm)\,w'(s_\pm)=\mp\infty.
\]
This matches the definition of an umbrella solution.
\end{proof}

We now show that, for fixed parameters, there exists a curve separating initial data that generate umbrella solutions from those that do not.

\begin{proposition}\label{proplambdaarriba}
Let $(a,\lambda,\delta,\gamma,\mu,c)$ be fixed and $\sigma\ge 0$.  
Then there exists a continuous curve
\[
\Lambda : 
\left[\frac{\sigma-c}{a},\,\frac{\sigma+c}{a}\right]
\longrightarrow \Theta_+,
\]
such that every initial datum above $\Lambda$ generates a trajectory of \eqref{eq-4} connecting $\pi_-$ and $\pi_+$.
\end{proposition}

\begin{proof}
We outline the construction.  
The possible sign changes of $w'$ correspond to up to three real roots $v_5<v_4<v_3$ of the equation $g(av-\sigma)-v=0$.  
Thus we partition $\Theta_+$ into the regions $\Theta_+^i$, $i=1,\dots,5$, determined by these values.  
Using the signs of $w'$ and the fact that $v'<0$ in $\Theta_+$, one checks that certain regions are positively invariant, while others are negatively invariant.

\medskip
\noindent\textbf{Step 1: Construction of one connecting orbit.}
Choosing sufficiently high initial data in $\Theta_+^1$, we use positive/negative invariance of the regions $\Theta_+^i$ and Proposition~\ref{propdesingular} to construct an orbit first reaching $v=v_5$, then $v=v_3$, and finally $\pi_+$.

\medskip
\noindent\textbf{Step 2: Monotone family of orbits.}
Fix $v_0$ and vary $w_0$ in the interval bounded above by the previously constructed orbit.  
Solutions depend monotonically on the initial height $w_0$, and the set of $w_0$ leading to connecting trajectories has a nonempty infimum.

\medskip
\noindent\textbf{Step 3: Definition and continuity of $\Lambda$.}
For each $v$ in the admissible interval, define $\Lambda(v)$ as the infimum height of the connecting orbits that pass through $v$.  
Using the monotonicity and compactness of the $v$-intervals comprising the regions $\Theta_+^i$, one proves that $\Lambda$ is continuous.

Thus $\Lambda$ is a continuous boundary separating initial data leading to umbrella solutions from initial data whose trajectories fail to connect $\pi_-$ and $\pi_+$.
\end{proof}

\begin{figure}\label{figPruebaConstruccion}
	\begin{tikzpicture}[scale=2]
		\draw[->] (-1,0)-- (3,0) node[right]{$v$} ;
		\draw[->] (1.5,0)--(1.5,2.5) node[above]{$w$};
		\draw[-,color=red, thick] (0,0) ..controls (1,1.5)..(2,0) ;
		\draw[dashed,color=red] (0,2.5) -- (0,0);
		\draw[ ] (0,0) node[below]{$P_1$} ;
		\draw[ ] (2,0) node[below]{$P_2$} ;
		
		\draw[ thick,color=red] (0.4,2.5) -- (0.4,0);
		\draw[] (0.4,0) node[below]{$v_5$};
		\draw[ thick,color=red] (1.2,2.5) -- (1.2,0);
		\draw[] (1.2,0) node[below]{$v_4$};
		\draw[ thick,color=red] (2.2,2.5) -- (2.2,0);
		\draw[] (2.2,0) node[below]{$v_3$};
		\draw[ thick,color=blue] (-0.8,2.5)node[right]{$\pi_-$} -- (-0.8,0);
		\draw[] (-0.8,0) node[below]{$\frac{\sigma-c}{a}$};
		\draw[ thick,color=blue] (2.8,2.5)node[right]{$\pi_+$} -- (2.8,0);
		\draw[] (2.8,0) node[below]{$\frac{\sigma+c}{a}$};
		
		\draw[fill=black] (2,0) circle(.03);
		\draw[fill=black] (0,0) circle(.03);
		\draw[fill=black] (0.4,0.59) circle(.03);
		\draw[fill=black] (1.2,1.05) circle(.03);
		
		\draw[-,color=black, thick] (-0.8,0.1) ..controls (-0.8,0.4) and (-0.2,0.2)..(0.4,0.8) ;
		\draw[-,color=black, thick] (0.4,1.7) ..controls (1.2,1.1) and (1.2,1.1)..(2.2,1.75) ;
		\draw[-,color=black, thick] (2.8,1.8) ..controls (2.7,1.9) and (2.4,2)..(2.2,2) ;
		\draw[-,color=green, thick] (2.2,2) ..controls (0,1.2) and (1.1,3)..(-0.8,1.2) ;
		
		
					\draw [->, thick] (0.1,1.1)--(-0.1,1) ;
					\draw [->, thick] (1.85,1.1)--(1.7,1) ;
					\draw [->, thick] (1,1.5)--(0.85,1.6) ;
					
					\draw [->, thick] (2.65,1.5)--(2.45,1.6) ;
					
					\draw [->, thick] (0.2,0.3)--(0.2,0.1);
					\draw [->, thick] (0.8,1.05)--(0.8,1.2);
					\draw [->, thick] (1.7,0.45)--(1.7,0.25);
		%
		%
	\end{tikzpicture}
	\caption{Esquema de demostración}
\end{figure}

Apart from the solutions described above, additional umbrella-type solutions may arise when studying orbits entirely contained in the region $\Theta_-$. 
The following result provides the corresponding separating curve.

\begin{proposition}\label{proplambdaabajo}
Let $(a,\lambda,\delta,\gamma,\mu,c)$ be fixed and $\sigma\ge 0$.  
Assume that the critical points $P_1$ and $P_2$ do not belong to the strip $\Theta$.  
Then there exists a continuous curve
\[
\Lambda:\ \left[\frac{\sigma-c}{a},\,\frac{\sigma+c}{a}\right]\longrightarrow \Theta_+,
\]
such that, for every initial datum $(v_0,w_0)$ lying below $\Lambda$, the system \eqref{eq-4} admits an orbit contained in $\Theta_-$ that connects the two singular lines $\pi_-$ and $\pi_+$.
\end{proposition}

\begin{proof}
Our aim is to construct connecting orbits lying entirely in $\Theta_-$.  
To make this possible, the equilibrium points
\[
P_1=(v_1,0), \qquad P_2=(v_2,0)
\]
must not lie inside the segment
\[
\Bigl[\tfrac{\sigma-c}{a},\,\tfrac{\sigma+c}{a}\Bigr]\times\{0\}.
\]
If either $P_1$ or $P_2$ belonged to this interval, any orbit in $\Theta_-$ would be attracted to that equilibrium, preventing a connection between $\pi_-$ and $\pi_+$.

Under this assumption, the construction follows the same strategy as in Proposition~\ref{proplambdaarriba}.  
The only difference is that, since $v'>0$ in $\Theta_-$, all invariant-set and comparison arguments must be performed from right to left, rather than left to right as in $\Theta_+$.  
This yields a continuous curve $\Lambda$ such that every initial condition below it generates a connecting orbit inside $\Theta_-$.
\end{proof}

With all these ingredients, we can now prove Theorem~\ref{teorema2}.

\begin{proof}[Proof of Theorem~\ref{teorema2}]
We prove the existence of the upper separating curve $\overline{\Lambda}$; the construction of the lower curve $\underline{\Lambda}$ is entirely analogous.

For every fixed $\sigma\ge0$, Proposition~\ref{proplambdaarriba} yields a continuous curve
\[
\Lambda_\sigma:\ 
\Bigl[\tfrac{\sigma-c}{a},\,\tfrac{\sigma+c}{a}\Bigr]
\longrightarrow \Theta_+, 
\qquad 
v\mapsto (v,\Lambda_\sigma(v)),
\]
such that all initial data lying above $\Lambda_\sigma$ generate an orbit connecting $\pi_-$ to $\pi_+$ while remaining in $\Theta_+$.  
Thus we obtain a family $\{\Lambda_\sigma\}_{\sigma\ge 0}$ indexed by the parameter~$\sigma$.

\textbf{Continuous dependence on the parameter.}  
Fix $\sigma_1\ge0$.  
A perturbation to $\sigma=\sigma_1+\varepsilon$ produces one of the following situations:

\begin{itemize}
\item[(i)] \emph{No qualitative change occurs.}  
The strip $\Theta$ shifts horizontally with $\sigma$, but the dynamical system undergoes no bifurcation: no equilibria are created, destroyed, or collide with singular boundaries.  
By performing a simple horizontal reparameterisation, the systems for $\sigma_1$ and $\sigma_1+\varepsilon$ may be written on the same domain, where the flows depend continuously on $\sigma$.  
Hence $\Lambda_\sigma$ depends continuously on $\sigma$ near $\sigma_1$.

\item[(ii)] \emph{A qualitative change occurs.}  
A local bifurcation takes place: an equilibrium is created or disappears, or two equilibria collide, or an equilibrium reaches a singular boundary.  
Such events occur only at isolated parameter values.  
Therefore, on any compact interval of $\sigma$, only finitely many values can produce bifurcations, and $\sigma\mapsto\Lambda_\sigma$ is continuous except at finitely many points.
\end{itemize}

\textbf{Domain of admissible $\sigma$ for fixed $v$.}  
A point $(v,w)$ belongs to $\Theta(\sigma)$ if and only if
\[
v\in\left(\tfrac{\sigma-c}{a},\,\tfrac{\sigma+c}{a}\right)
\quad\Longleftrightarrow\quad
\sigma\in\bigl[\max\{0,av-c\},\, av+c\bigr].
\]
Thus, for a fixed $v>-c/a$, only the compact interval 
\[
\Sigma(v):=\bigl[\max\{0,av-c\},\, av+c\bigr]
\]
is relevant when comparing the curves~$\Lambda_\sigma$.

\textbf{Construction of the upper separating curve.}  
For each fixed $v>-c/a$, define
\[
\overline{\Lambda}(v)
:=\sup_{\sigma\in\Sigma(v)} \Lambda_\sigma(v).
\]
This supremum is finite because all curves $\Lambda_\sigma$ lie inside $\Theta_+$.  
Since $\Lambda_\sigma(v)$ is continuous in $\sigma$ except at finitely many values, 
the supremum over a compact interval $\Sigma(v)$ defines a well-behaved function of~$v$.  
Moreover, as $v$ varies, discontinuities can only occur at those finitely many $v$ where the maximizing parameter jumps from one side of a bifurcation value to the other.  
Thus $\overline{\Lambda}(v)$ is continuous except at finitely many points of its domain.

\textbf{Conclusion.}  
By construction,
\[
\Lambda_\sigma(v)\le \overline{\Lambda}(v)
\qquad\text{for all admissible }\sigma\in\Sigma(v).
\]
Hence any initial condition $(v,w)$ with $w>\overline{\Lambda}(v)$ lies above \emph{every} curve $\Lambda_\sigma$, and therefore, for every admissible $\sigma$, the corresponding orbit connects $\pi_-$ to $\pi_+$.

This completes the proof.
\end{proof}

\end{document}